\documentclass[11pt,reqno]{elsarticle}

\makeatletter
\def\ps@pprintTitle{%
 \let\@oddhead\@empty
 \let\@evenhead\@empty
 \def\@oddfoot{}%
 \let\@evenfoot\@oddfoot}
\makeatother

\usepackage{amsmath,amsthm,amscd,amsfonts,amssymb,mathtools}
\usepackage{mathabx,breqn,multirow,comment}
\usepackage{booktabs,cool}
\usepackage{parskip}
\usepackage[top=1.1in, bottom=1.1in, left=1in, right=1in]{geometry}

\newtheorem{remark}{Remark}

\providecommand{\e}[1]{\ensuremath{\times 10^{#1}}}

\makeatletter
\g@addto@macro\normalsize{%
  \setlength\abovedisplayskip{.4em}
  \setlength\belowdisplayskip{.4em}
  \setlength\abovedisplayshortskip{.4em}
  \setlength\belowdisplayshortskip{.4em}
}

\begin{document}

\begin{frontmatter}
\begin{abstract}
We describe a method for the rapid numerical evaluation of the Bessel 
functions of  the first and second kinds of nonnegative real orders and 
positive arguments.    Our algorithm makes use of the well-known observation
that although the Bessel functions themselves are expensive to represent
via piecewise polynomial expansions, the logarithms of certain solutions
of Bessel's equation are not.  We exploit this observation by numerically
precomputing the logarithms of carefully chosen Bessel functions  and representing them
with piecewise bivariate Chebyshev expansions.  
We supplement these precomputed expansions with two asymptotic expansions,
one for large orders and extremely small arguments and the other for large orders
and extremely large arguments, and with series expansions  for  small orders and small arguments.  
Our scheme is able to evaluate Bessel functions
of orders between $0$ and $1\sep,000\sep,000\sep,000$
at essentially any positive real argument. 
In that regime, it is competitive with existing methods for the rapid
evaluation of Bessel functions and has at least three advantages
over them.   First, our approach is quite general
and can be readily applied to many other special functions which satisfy second
order ordinary differential equations.
Second, by calculating the logarithms of the Bessel functions rather
than the Bessel functions themselves, we avoid many issues which arise from numerical overflow
and underflow.  Third, in the oscillatory regime, our algorithm calculates the values of a nonoscillatory
phase function for Bessel's differential equation and its derivative.  These quantities
are useful for computing the zeros of Bessel functions, as well as  for rapidly applying
the Fourier-Bessel transform.
The results of extensive numerical experiments demonstrating the efficacy of our algorithm
are presented.  A Fortran package which includes our code for evaluating
the Bessel functions as well  as our code for all of the numerical experiments 
described here is publically available.

\end{abstract}

\begin{keyword}
special functions \sep
fast algorithms \sep
nonoscillatory phase functions
\end{keyword}

\title
{
An algorithm for the rapid numerical evaluation of Bessel functions
of real orders and arguments
}

\author[jb]{James Bremer}
\ead{bremer@math.ucdavis.edu}

\address[jb]{Department of Mathematics, University of California, Davis}

\end{frontmatter}

\begin{section}{Introduction}

Here, we describe a numerical method for evaluating
the Bessel functions of the first and second kinds ---
$J_\nu$ and $Y_\nu$, respectively --- of nonnegative orders and positive arguments.  
In this regime, it is competitive with (and possess some advantages over) 
existing methods for the numerical evaluation
of the Bessel functions such as \cite{Amos} and \cite{Matviyenko}.

The purpose of this article, though, is not to argue that existing
schemes for the evaluation of Bessel functions
are inadequate or should be replaced with ours.
Instead, it is to point out  that there is an incredibly
straightforward  approach to their numerical evaluation  that applies to a large
class of special functions satisfying second order 
ordinary differential equation.  Our decision to focus in the first instance
on Bessel functions  stems in large part from the existence of satisfactory
numerical algorithms with which we can compare our approach.
The results of applying  this approach to other classes of special functions, such
as the associated Legendre functions and prolate spheroidal wave functions,
will be reported by the author at a later date.

It is well known that the scaled Bessel functions $J_\nu(t) \sqrt{t}$ and $Y_\nu(t) \sqrt{t}$ satisfy
the second order linear ordinary differential equation
\begin{equation}
y''(t) + \left( 1 - \frac{\nu^2 - \frac{1}{4}}{t^2} \right) y(t) = 0 \ \ \mbox{for all} \ \ 0 < t< \infty.
\label{introduction:besseleq}
\end{equation}
We will, by a slight abuse of terminology, refer to (\ref{introduction:besseleq})
as Bessel's differential equation.  When $0 \leq \nu \leq 1/2$, the coefficient of
$y$ in (\ref{introduction:besseleq}) is positive on the entire half-line
$(0,\infty)$, whereas it is negative on the interval
\begin{equation}
\left(0, \sqrt{\nu^2-\frac{1}{4}}\right)
\label{introduction:interval1}
\end{equation}
and positive on
\begin{equation}
\left(\sqrt{\nu^2-\frac{1}{4}},\infty\right)
\label{introduction:interval2}
\end{equation}
when $\nu > 1/2$.    
It follows from  standard WKB estimates  (see, for instance, \cite{Fedoryuk})
that  solutions of (\ref{introduction:besseleq})
 behave roughly like  increasing or decreasing exponentials on (\ref{introduction:interval1}) 
and are oscillatory on (\ref{introduction:interval2}). 
We will refer to the subset
\begin{equation}
\mathcal{O} = 
\left\{
(\nu,t) :
0 \leq \nu \leq \frac{1}{2}\ \ \mbox{and} \ t > 0
\right\}
\bigcup
\left\{
(\nu,t) : \nu > \frac{1}{2}\ \ \mbox{and}\ 
t \geq \sqrt{\nu^2-\frac{1}{4}}
\right\}
\label{introduction:oscillatory}
\end{equation}
of $\mathbb{R} \times \mathbb{R}$ as the oscillatory region and the subset
\begin{equation}
\mathcal{N} = 
\left\{
(\nu,t) :
\nu > \frac{1}{2}  \ \ \mbox{and}\ \ 
0 < t < \sqrt{\nu^2-\frac{1}{4}}
\right\}
\label{introduction:nonoscillatory}
\end{equation}
of $\mathbb{R}\times \mathbb{R}$ as the nonoscillatory region.
%


When $\nu$ is large, we cannot expect to represent Bessel functions
efficiently using polynomial expansions in $\mathcal{N}$ because, in this event,
they behave like rapidly increasing or decreasing exponentials.
Similarly, we cannot expect to represent Bessel functions efficiently
with polynomial expansions on any substantial subset of $\mathcal{O}$ 
since they oscillate there.
Despite this, the logarithms of  Bessel functions can be 
represented efficiently via polynomial expansions on $\mathcal{N}$.
Moreover, there is a \emph{carefully selected} solution of Bessel's differential equation
whose logarithm can be represented efficiently  via polynomial expansions on a substantial
subset of 
oscillatory region $\mathcal{O}$.  This latter observation is related to the well-known
fact that Bessel's differential equation admits a nonoscillatory phase
function (see, for instance, Section~13.75 of \cite{Watson} or
\cite{Heitman-Bremer-Rokhlin-Vioreanu}).

Many special functions of interest share this property of Bessel functions,
at least in an asymptotic sense \cite{Miller,DLMF}.  However, 
the sheer effectiveness with which nonoscillatory phase functions
can represent solutions of the general equation
\begin{equation}
y''(t) + \lambda^2 q(t) y(t) = 0 \ \ \mbox{for all}\ \ a < t <b
\label{introduction:second_order}
\end{equation}
in  which the coefficient $q$ is smooth and positive 
appears to have been overlooked.  Indeed, under mild conditions on $q$, 
it is shown in  \cite{Bremer-Rokhlin} that there exist a positive real number $\mu$,
 a nonoscillatory function $\alpha$  and a basis
of solutions $\{u,v\}$ of (\ref{introduction:second_order}) such that
\begin{equation}
u(t) = \frac{\cos\left(\alpha(t)\right)}{\sqrt{\alpha'(t)}} + O\left(\exp(-\mu\lambda)\right)
\end{equation}
and
\begin{equation}
v(t) = \frac{\sin\left(\alpha(t)\right)}{\sqrt{\alpha'(t)}} + O\left(\exp(-\mu\lambda)\right).
\end{equation}
The constant $\mu$ is a measure of the extent to which $q$ oscillates, with
larger values of $\mu$ indicating greater smoothness on the part of $q$.
The function $\alpha$ is nonoscillatory in the sense that it can be represented
using various series expansions the number of terms in which do not vary
with $\lambda$.  That is, $O(\exp(-\mu\lambda))$ accuracy is obtained using
an $O(1)$-term expansion.
The  results of \cite{Bremer-Rokhlin} are akin to standard results on
WKB approximation in that they apply to the more general case in which
$q$ varies with the parameter $\lambda$ assuming only that
$q$ satisfies certain innocuous hypotheses independent of 
$\lambda$.  An effective numerical algorithm for the computation
of nonoscillatory phase functions for equations of the form (\ref{introduction:second_order})
is described in \cite{BremerKummer}, although we will not need it here
since an effective asymptotic expansion of a nonoscillatory phase function
for Bessel's differential equation is available.
The algorithm of \cite{BremerKummer} is, however, of importance in generalizing
these results to cases in which such expansions are not available.

The algorithm of this paper 
operates by numerically precomputing the logarithms of certain
solutions of Bessel's differential equation.  We represent these
function via piecewise bivariate Chebyshev expansions,  the coefficients
of which are arranged in a table.  The table is, of course, stored on the disk
and loaded into memory when needed so it only needs to be computed once.
We supplement these precomputed expansions with asymptotic and series expansions in order
to evaluate  $J_\nu$ and $Y_\nu$ for all nonnegative real orders $\nu$
and  positive real arguments.  We note, though, that
in cases in which such expansions are not available,
the range of the parameter and argument covered by the precomputed expansions is 
 sufficient for most purposes and could be extended if needed.
The size of the precomputed table used in the experiments of
this paper is roughly $1.3$ megabytes.

The remainder of this paper is structured as follows.  In Section~\ref{section:preliminaries},
we review certain mathematical facts and numerical procedures which are used in the
rest of this  article.  Section~\ref{section:solver} details the operation
of a solver for nonlinear differential equations
which is used by the algorithm of Section~\ref{section:phase}
for the rapid  solution of Bessel's differential equation
(\ref{introduction:besseleq}) in the case in which the parameter $\nu$ is fixed.  
This procedure is, in turn, a component of the scheme for the construction
of the precomputed table which we use to evaluate Bessel functions.
That scheme is described in  Section~\ref{section:expansions}.
Section~\ref{section:numerics} details our algorithm for the numerical
evaluation of Bessel functions using this table and certain asymptotic and series
expansions.  Section~\ref{section:experiments} describes extensive numerical experiments 
performed  in order to verify the efficacy of the algorithm of 
Section~\ref{section:numerics}.  We conclude with a few remarks regarding this contents
of this article and possible directions for future work in Section~\ref{section:conclusion}.
\end{section}

\begin{section}{Mathematical and Numerical Preliminaries}


\begin{subsection}{The condition number of the evaluation of a function}

The condition number of the evaluation of a differentiable
function $f:\mathbb{R} \to\mathbb{R}$
at the point $x$  is commonly defined to be 
\begin{equation}
\kappa_{f}(x) = \left| \frac{x f'(x) }{f(x)} \right|
\label{preliminaries:condition:1}
\end{equation}
(see, for instance, Section~1.6 of \cite{Higham}). 
This quantity measures the ratio of the magnitude of 
the relative change in $f(x)$ induced by a small change
in the argument $x$ to the magnitude of the relative change in $x$
in the sense that 
\begin{equation}
\left| \frac{f(x+\delta) - f(x)}{f(x)} \right|
\approx  \kappa_{f}(x)\ \left|\frac{\delta}{x}\right|
\end{equation}
for small $\delta$.  Since almost all quantities which arise in the course
of numerical calculations  are subject to perturbations
with relative magnitudes on the order of machine epsilon,
we consider 
\begin{equation}
\kappa_f(x) \epsilon_0,
\end{equation}
where $\epsilon_0$ denotes machine epsilon, to be a rough estimate of the 
relative accuracy one should expect when evaluating $f(x)$ numerically
(in fact, it tends to be a  slightly pessimistic estimate).
In the rest of this paper, we take $\epsilon_0$ to be 
\begin{equation}
\epsilon_0 = 2^{-52} \approx 2.22044604925031\e{-16}.
\end{equation}
It is immediately clear from  (\ref{preliminaries:condition:1}) that
when $f'(x_0) x_0 \neq 0$ and $f(x_0) = 0$, $\kappa_{f}(x)$ diverges to $\infty$ as  $x \to x_0$.
One consequence of this is that there is often a significant loss
of relative accuracy when  $f(x)$ is evaluated near one of its roots.
In order to avoid this issue, we will arrange for all of the quantities we calculate
 to be bounded away from $0$.

\label{preliminaries:condition}
\end{subsection}


\begin{subsection}{Series expansions of the Bessel functions}

For complex-valued $\nu$ and $x>0$, the Bessel function of the first kind
of order $\nu$ is given by
\begin{equation}
J_\nu(x) = \sum_{j=0}^\infty \frac{(-1)^j}{\Gamma(j+1) \Gamma(j+\nu+1)} 
\left(\frac{x}{2}\right)^{2j+\nu}.
\label{preliminaries:series:1}
\end{equation}
Here, we use the convention that 
\begin{equation}
\frac{1}{\Gamma(j)} = 0
\end{equation}
whenever $j$ is a negative integer.  Among other things, this ensures
 that (\ref{preliminaries:series:1}) is still sensible when $\nu$ is a negative integer.
When $x > 0$ and $\nu$ is not an integer,
the Bessel function of the second kind  of order $\nu$ is given by
\begin{equation}
Y_\nu(x) =  \frac{  \cos(\nu\pi) J_\nu(x) - J_{-\nu}(x) } {\sin(\nu\pi)}.
\label{preliminaries:series:2}
\end{equation}
For integer values of $\nu$, (\ref{preliminaries:series:2}) loses its meaning;
however, taking the limit of $Y_\nu(x)$ as $\nu \to n \in \mathbb{Z}$ yields
\begin{equation}
\begin{aligned}
Y_n(x) = \frac{2}{\pi} J_n(x) &\log \left(\frac{x}{2}\right)
-  \sum_{j=0}^{n-1} \frac{\Gamma(n-j-1)}{\Gamma(j+1)}  \left(\frac{x}{2}\right)^{2j-n} \\
-& \frac{1}{\pi}
\sum_{j=0}^\infty (-1)^j\ \frac{\psi(n+j+1) + \psi(j+1)}{\Gamma(j+1) \Gamma(n+j+1)}  
\left(\frac{x}{2}\right)^{n+2j},
\end{aligned}
\label{preliminaries:series:3}
\end{equation}
where $\psi$ is the logarithmic derivative of the gamma function.
A derivation of this formula can be found in Section~7.2.4 of \cite{HTFII}.

For the most part, when $\nu$ is of small in magnitude and $t$ is positive
and of small magnitude, the value of $J_\nu(t)$
can be computed in a numerically stable fashion by truncating the series
(\ref{preliminaries:series:1}).  In some cases, however, this
can lead to numerical underflow.  Accordingly, when $t \ll \nu$, we 
evaluate the logarithm of $J_\nu$ by truncating the series in
the expression
\begin{equation}
\log(J_\nu(t)) = 
-\log(\Gamma(\nu+1)) + \nu \log\left(\frac{t}{2}\right) + 
\log\left(\sum_{j=0}^\infty \frac{(-1)^j \Gamma(\nu+1)}{\Gamma(j+1) \Gamma(j+\nu+1)}
\left(\frac{x}{2}\right)^j \right)
\label{preliminaries:series:4}
\end{equation}
instead.

On the other hand, Formula~(\ref{preliminaries:series:2}) can lead to significant errors 
when it used to evaluate $Y_\nu$ numerically.
In particular,  when $\nu$ is close to, but still distinct from an integer,
the evaluation of $Y_\nu$ via (\ref{preliminaries:series:2})
results in significant round-off error due to numerical cancellation.
Since $Y_\nu$ is analytic as a function of $\nu$, this problem can be
obviated by evaluating $Y_\nu$ via interpolation with respect to the
order $\nu$.     Similarly, it is often more convenient
to compute $Y_\nu$ when $\nu$ is an integer using interpolation
than to do so via (\ref{preliminaries:series:3}).  Similar suggestions are 
made in \cite{Matviyenko}.

The naive use  of (\ref{preliminaries:series:2}) can also lead to numerical overflow
 when $\nu$ is not close to an integer.  In such cases we evaluate $\log(-Y_\nu(t))$ via 
\begin{equation}
\log(-Y_\nu(x)) =  
\log\left(J_\nu(x)\right) +
\log\left(
\frac{  -\cos(\nu\pi)  + 
\exp(\log(J_{-\nu}(x)) - \log(J_\nu(x)))
 } {\sin(\nu\pi)}
\right).
\label{preliminaries:series:5}
\end{equation}
We calculate the logarithms of $J_\nu$ appearing in (\ref{preliminaries:series:5}) using
(\ref{preliminaries:series:4}), of course.

\label{preliminaries:series}
\end{subsection}


\begin{subsection}{Debye's asymptotic expansion for small arguments}

The following form of Debye's asymptotic expansions can be found
in \cite{Matviyenko}.     For $x < \nu$ and $N$ a nonnegative integer,
\begin{equation}
J_\nu(x) = \frac{1}{1+ \theta_{N+1,1}(\nu,0)} 
\frac{\exp(-\eta)}{\sqrt{2\pi} (\nu^2-x^2)^{\frac{1}{4}}} \times
\left(
\sum_{j=0}^N \frac{u_j(p)}{\nu^j}
+ \theta_{N+1,1}(\nu,p)
\right)
\label{preliminaries:debye:expansion1}
\end{equation}
and
\begin{equation}
Y_\nu(x) = - \sqrt{\frac{2}{\pi}}
\frac{\exp(\eta)}{(\nu^2-x^2)^{\frac{1}{4}}} \times
\left( 
\sum_{j=0}^N (-1)^j \frac{u_j(p)}{\nu^j}
+ \theta_{N+1,2}(\nu,p)
\right),
\label{preliminaries:debye:expansion2}
\end{equation}
where 
\begin{equation}
\eta = \nu \log\left(
\frac{\nu}{x}  + 
\sqrt{ \left(\frac{\nu}{x}\right)^2 -1 } 
\right)
- \sqrt{\nu^2-x^2},
\end{equation}
\begin{equation}
p = \frac{\nu}{\sqrt{\nu^2-x^2}},
\end{equation}
$\theta_{N+1,1}$ and $\theta_{N+1,2}$ are error terms, 
and $u_0$, $u_1$,  $\ldots$ are the polynomials defined via
\begin{equation}
u_0(t) = 1,
\end{equation}
and the recurrence relation
\begin{equation}
u_{n+1}(t) = \frac{1}{2} \left(t^2 -t^4\right) \frac{du_n(t)}{dt} +
\frac{1}{8} \int_0^t (1-5\tau^2)u_n(\tau)\ d\tau
\ \ \mbox{for all} \ \ n \geq 0.
\end{equation}
In  \cite{Matviyenko}, it is shown that there exist positive real constants
$C_1,C_2,\ldots$ such that
\begin{equation}
\max\left\{\left|\theta_{N+1,1}(\nu,p)\right|,
\left|\theta_{N+1,2}(\nu,p)\right|
\right\}
 \leq 2 \exp\left(\frac{2}{3 g^{\frac{3}{2}}}\right)
\frac{C_{N+1}}{g^{\frac{3}{2} (N+1)}},
\end{equation}
where
\begin{equation}
g = \frac{\nu-x}{\nu^{\frac{1}{3}}},
\label{preliminaries:debye:1}
\end{equation}
for all $N \geq 0$.  In other words, Debye's asymptotic expansions for small 
values of the parameter are uniform asymptotic expansions in inverse powers of the variable
(\ref{preliminaries:debye:1}).  See \cite{Matviyenko} for a further discussion
of the implications of this observation.

The naive use of (\ref{preliminaries:debye:expansion1}) and (\ref{preliminaries:debye:expansion2})
when $t \ll \nu$  often results in numerical underflow and overflow.  In order to avoid such problems,
in this regime
we evaluate the logarithms of the  Bessel functions  via the approximations
\begin{equation}
\log\left(J_\nu(x)\right) \approx
-\eta - \frac{1}{4} \log(\nu^2-x^2)
+ \log\left(\frac{1}{\sqrt{2\pi}}
\sum_{j=0}^N \frac{u_j(p)}{\nu^j}
\right)
\label{preliminaries:debye:expansion3}
\end{equation}
and
\begin{equation}
\log(-Y_\nu(x)) 
\approx
\eta  - \frac{1}{4} \log(\nu^2-x^2)
+ \log \left(
 \sqrt{\frac{2}{\pi}}
\sum_{j=0}^N (-1)^j \frac{u_j(p)}{\nu^j}
\right)
\label{preliminaries:debye:expansion4}
\end{equation}
rather than evaluate the Bessel functions themselves.

\label{preliminaries:debye}
\end{subsection}


\begin{subsection}{The Riccati equation, Kummer's equation and phase functions}

If $y = \exp(r(t))$ satisfies
\begin{equation}
y''(t) + q(t) y(t) = 0 \ \ \mbox{for all}\ \ t \in I,
\label{preliminaries:riccati:1}
\end{equation}
where $I \subset \mathbb{R} $ is an open interval,
then a straightforward computation shows that
\begin{equation}
r''(t) + (r(t))^2 + q(t) = 0\ \  \mbox{for all}\ \ t\in I.
\label{preliminaries:riccati:2}
\end{equation}
 Equation~(\ref{preliminaries:riccati:2}) is known as the Riccati equation;
an extensive discussion of it can be found, for instance, in \cite{Hille}.
By assuming that $q$ is real-valued, and that
\begin{equation}
r(t) = \alpha(t) + i \beta(t)
\label{preliminaries:riccati:3}
\end{equation}
with $\alpha$ and $\beta$ real-valued, we obtain 
from (\ref{preliminaries:riccati:2}) the system of ordinary
differential equations
\begin{equation}
\left\{
\begin{aligned}
\beta''(t) + (\beta'(t))^2 - (\alpha'(t))^2 + q(t) &= 0 \\
\alpha''(t) + 2 \alpha'(t) \beta'(t) &= 0.
\end{aligned}
\right.
\label{preliminaries:riccati:4}
\end{equation}
If $\alpha'$ is nonzero, then the second of these equations readily implies that
\begin{equation}
\beta(t) = -\frac{1}{2} \log\left(\left|\alpha'(t)\right|\right).
\label{preliminaries:riccati:5}
\end{equation}
Inserting (\ref{preliminaries:riccati:5}) into the first equation
in (\ref{preliminaries:riccati:4}) yields
\begin{equation}
 q(t) 
- (\alpha'(t))^2
- \frac{1}{2}\left(\frac{\alpha'''(t)}{\alpha'(t)}\right)
+ \frac{3}{4}
\left(\frac{\alpha''(t)}{\alpha'(t)}\right)^2 = 0
\label{preliminaries:riccati:kummer}.
\end{equation}
We will refer to (\ref{preliminaries:riccati:kummer}) as Kummer's equation
after E.~E. Kummer who studied it in \cite{Kummer}.
We conclude  that if the derivative of the function $\alpha$ is nonzero
and  satisfies (\ref{preliminaries:riccati:kummer}), then 
\begin{equation}
u(t) = \frac{\cos\left(\alpha(t)\right)}{\sqrt{\left|\alpha'(t)\right|}}
\label{preliminaries:riccati:u}
\end{equation}
and
\begin{equation}
v(t) = \frac{\sin\left(\alpha(t)\right)}{\sqrt{\left|\alpha'(t)\right|}}
\label{preliminaries:riccati:v}
\end{equation}
are   solutions of the differential equation
(\ref{preliminaries:riccati:1}).     A straightforward computation
shows that the Wronskian of $\{u,v\}$ is $1$, so that 
they form a basis in the space of solutions of this differential equation.
In this event, the function $\alpha$
is said to be  a phase function for (\ref{preliminaries:riccati:1}).

Suppose, on the other hand, that $\tilde{u}$ and $\tilde{v}$ are real-valued 
 solutions of (\ref{preliminaries:riccati:1}), that the Wronskian
of $\{\tilde{u},\tilde{v}\}$ is $1$, and that $\alpha$ is a smooth function such that
\begin{equation}
\alpha'(t) = \frac{1}{\tilde{u}(t)^2 + \tilde{v}(t)^2}.
\label{preliminaries:riccati:alphap}
\end{equation}
Since $\tilde{u}$ and $\tilde{v}$ cannot simultaneously vanish on $I$, 
$\alpha'$ is necessarily positive there.  A tedious, but straightforward computation
shows that (\ref{preliminaries:riccati:alphap}) satisfies Kummer's equation
(\ref{preliminaries:riccati:kummer}) so that $\alpha$ 
is a phase function for (\ref{preliminaries:riccati:1}) and 
the functions $u$, $v$ defined via
(\ref{preliminaries:riccati:u}) and (\ref{preliminaries:riccati:v})
form a basis in its space of solutions.
  We note, though, that
since (\ref{preliminaries:riccati:alphap}) only determines $\alpha$ up to a constant,
$u$ need not coincide with $\tilde{u}$ and $v$ need not coincide with $\tilde{v}$.

\label{preliminaries:riccati}
\end{subsection}

\begin{subsection}{A nonoscillatory phase function for Bessel's equation}

In the case of the solutions
\begin{equation}
u_\nu(t) = \sqrt{\frac{\pi t}{2}  } J_\nu(t)
\label{preliminaries:nonoscillatory:1}
\end{equation}
and
\begin{equation}
v_\nu(t) = \sqrt{\frac{\pi t}{2}  } Y_\nu(t)
\label{preliminaries:nonoscillatory:2}
\end{equation}
of Bessel's differential equation,   (\ref{preliminaries:riccati:alphap}) becomes
\begin{equation}
\alpha_\nu'(t) =  \frac{2}{\pi t}\ \frac{1}{J_\nu^2(t) + Y_\nu^2(t)}.
\label{preliminaries:nonoscillatory:alphap}
\end{equation}
Note that the  Wronskian of the pair  $\{u_\nu,v_\nu\}$ is $1$ on the interval $(0,\infty)$  
(see, for instance, Formula~(28) in Section~7.11  of \cite{HTFII}).
We define a phase function $\alpha_\nu$ for (\ref{introduction:besseleq}) via the formula
\begin{equation}
\alpha_\nu(t)  = C + \int_0^t \alpha_\nu'(s)\ ds
\end{equation}
with the constant $C$ to be set so that
\begin{equation}
\frac{\cos(\alpha_\nu(t))}{\sqrt{\alpha_\nu'(t)}} = u_\nu(t)
\label{preliminaries:nonoscillatory:3}
\end{equation}
and
\begin{equation}
\frac{\sin(\alpha_\nu(t))}{\sqrt{\alpha_\nu'(t)}} = v_\nu(t).
\label{preliminaries:nonoscillatory:4}
\end{equation}
  From (\ref{preliminaries:nonoscillatory:alphap})
and the series expansions for $J_\nu$ and $Y_\nu$  
appearing in Section~\ref{preliminaries:series},  we see that 
\begin{equation}
\lim_{t \to 0^+} \sqrt{\alpha_\nu'(t) }\  u_\nu(t) = 0 
\end{equation}
while
\begin{equation}
\lim_{t \to 0^+} \sqrt{\alpha_\nu'(t) }\  v_\nu(t) = -1.
\end{equation}
It follows that in order for
\begin{equation}
\lim_{t\to 0^+}\ \left( \frac{\cos(\alpha_\nu(t)) }{\sqrt{\alpha_\nu'(t)}} - u_\nu(t) \right) = 0 =
\lim_{t\to 0^+}\ \left( \frac{\sin(\alpha_\nu(t)) }{\sqrt{\alpha_\nu'(t)}} - v_\nu(t)\right) 
\end{equation}
to hold, we must have
\begin{equation}
\cos(C)=\cos(\alpha_\nu(0)) = 0
\end{equation}
and
\begin{equation}
\sin(C)=\sin(\alpha_\nu(0)) = -1.
\end{equation}
We conclude that by taking $C = -\pi/2$ --- so that
\begin{equation}
\alpha_\nu(t)  = -\frac{\pi}{2} + \int_0^t \alpha_\nu'(s)\ ds
\label{preliminaries:nonoscillatory:alpha}
\end{equation}
--- we ensure that (\ref{preliminaries:nonoscillatory:3}) and (\ref{preliminaries:nonoscillatory:4})
are satisfied.

\begin{figure}[b!]
\begin{center}
\includegraphics[width=.40\textwidth]{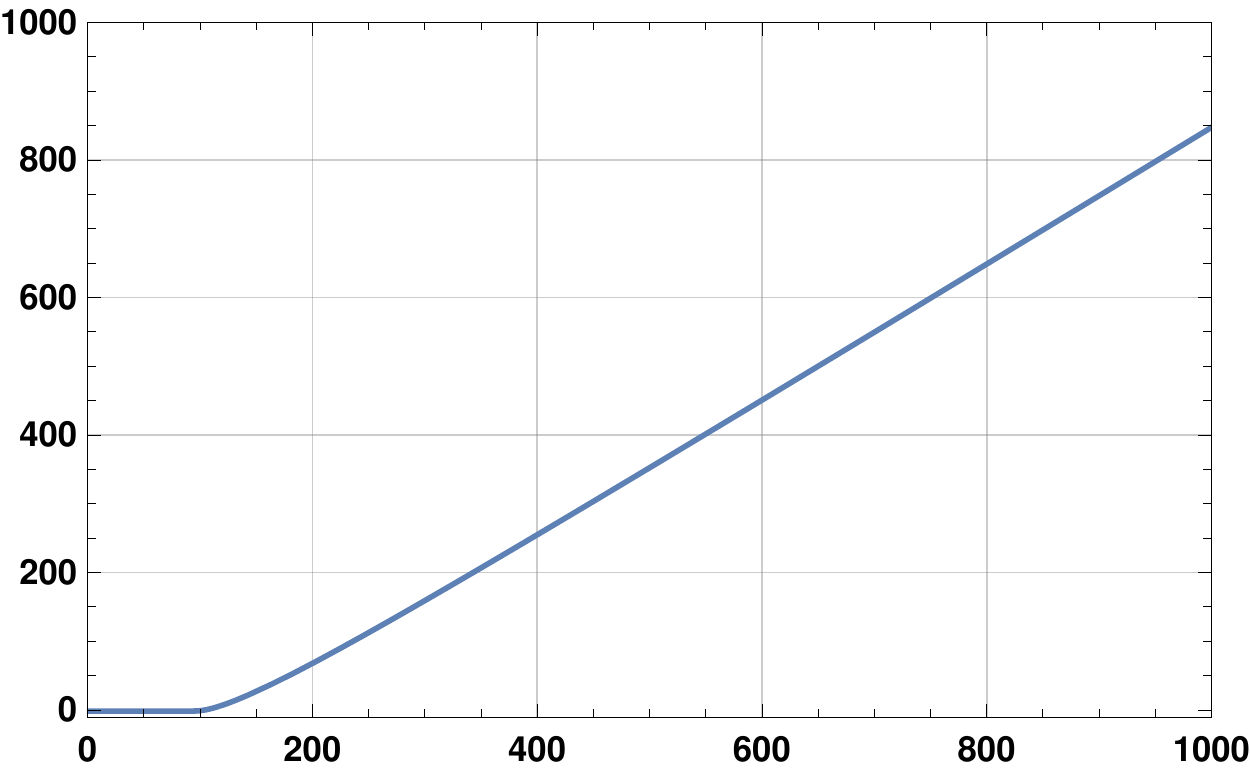}
\hfil
\includegraphics[width=.40\textwidth]{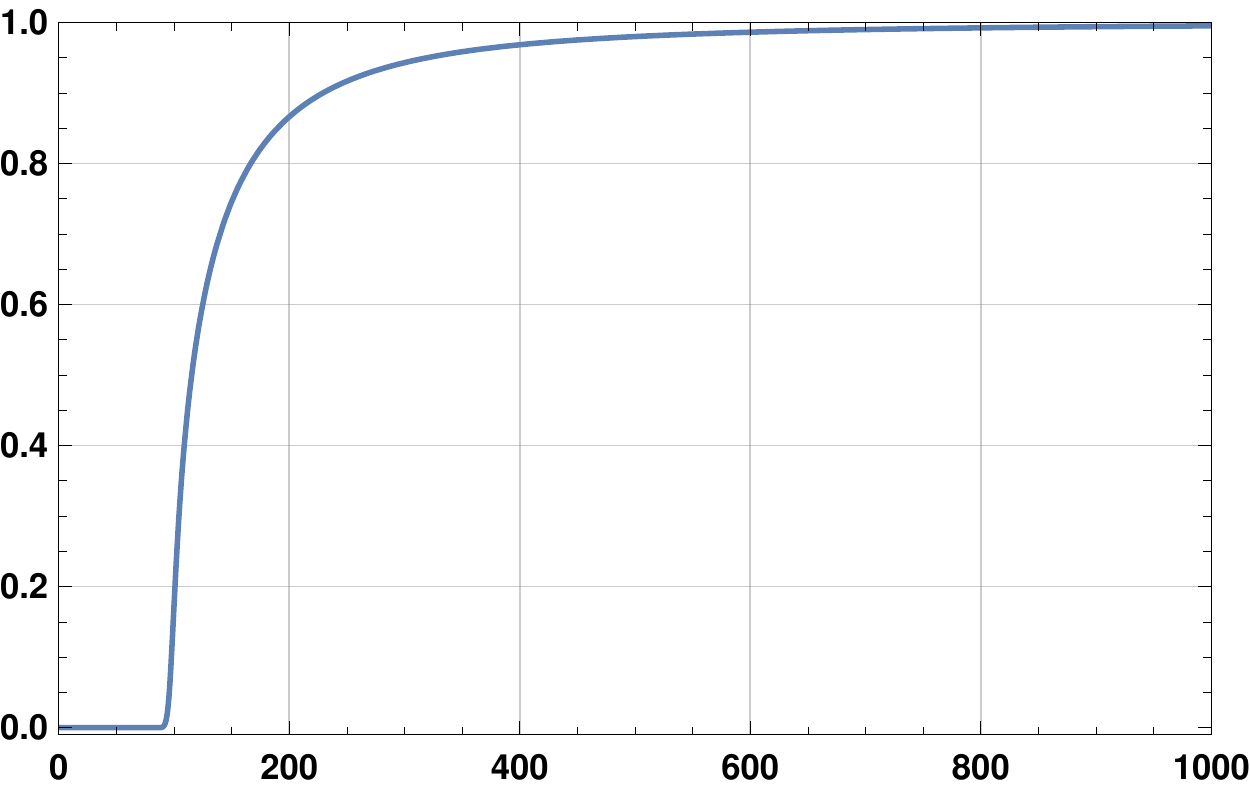}
\end{center}
\caption{
On the left is a plot of the phase function $\alpha_\nu(t)$ defined
via (\ref{preliminaries:nonoscillatory:alpha}) when $\nu=100$,
and on the right is a plot of its derivative with respect to $t$.
}
\label{preliminaries:figure1}
\end{figure}

Now we denote by $M_\nu$ the function appearing in denominator of the second
factor in (\ref{preliminaries:nonoscillatory:alphap}); that is, $M_\nu$ is defined via
\begin{equation}
M_\nu(t) = (J_\nu(t))^2 + (Y_\nu(t))^2.
\end{equation}
A cursory examination of Nicholson's formula
\begin{equation}
M_\nu(t) = \frac{8}{\pi^2} \int_0^\infty K_0(2t\sinh(s)) \cosh(2\nu s)\ ds
\ \ \mbox{for all}\ \  t >0,
\label{preliminaries:nonoscillatory:nicholson}
\end{equation}
a derivation of which can be found in Section~13.73 of \cite{Watson},
shows that $M_\nu$, and hence also $\alpha_\nu'$ and $\alpha_\nu$,
are nonoscillatory as  functions of $t$.
Figure~\ref{preliminaries:figure1} shows plots of the nonoscillatory functions
$\alpha_\nu$ and $\alpha_\nu'$ when $\nu = 100$.
We note that while $\alpha_\nu$ and $\alpha_\nu'$ can be represented
efficiently via polynomial expansions in the oscillatory regime, it is
clear from these plots that that is not the case in the nonoscillatory
regime.  
Moreover, the asymptotic expansion
\begin{equation}
M_\nu(t) \sim
\frac{2}{\pi t} 
\sum_{n=0}^\infty \frac{\Gamma\left(n+\frac{1}{2}\right)}{\sqrt{\pi}\ \Gamma(n+1)}
\frac{\Gamma\left(\nu + n + \frac{1}{2} \right)}{\Gamma\left(\nu - n +  \frac{1}{2} \right)}
\frac{1}{t^{2n}}
\ \ \mbox{as}\ \ t\to\infty
\label{preliminaries:nonoscillatory:5}
\end{equation}
can be derived easily from (\ref{preliminaries:nonoscillatory:nicholson}) 
(see Section~13.75 of \cite{Watson}).    The first few terms are
\begin{equation}
M_\nu \sim
\frac{2}{\pi t} \left( 1 + \frac{1}{2} \frac{\mu-1}{(2t)^2}  
+ \frac{1}{2}\cdot \frac{3}{4} \frac{(\mu-1)(\mu-9)}{(2t)^4} 
+ \frac{1}{2}\cdot \frac{3}{4} \cdot \frac{5}{6} \frac{(\mu-1)(\mu-9)(\mu-25)}{(2t)^6} + 
\cdots
\right),
\label{preliminaries:nonoscillatory:6}
\end{equation}
where $\mu = 4\nu^2$.    

From (\ref{preliminaries:nonoscillatory:5}), we can
derive an asymptotic expansion of (\ref{preliminaries:nonoscillatory:alphap}).
Indeed, if we denote the $n^{th}$ coefficient in the sum appearing in
(\ref{preliminaries:nonoscillatory:5})  by $r_n$, 
then the coefficients $s_0,s_1,\ldots$ in the asymptotic expansion
\begin{equation}
\alpha_\nu'(t) \sim \sum_{n=0}^\infty \frac{s_n}{t^{2n}}
\ \ \mbox{as}\ \ t \to\infty
\label{preliminaries:nonoscillatory:alphap_asym}
\end{equation}
are given by
\begin{equation}
s_0 = 1 \ \ \mbox{and} \ \ 
s_{n} = -\sum_{j=1}^{n} s_{n-j} r_{j} \ \ \mbox{for all}\ \ n \geq 1.
\label{preliminaries:nonoscillatory:7}
\end{equation}
The proof of this is an exercise in elementary calculus, and can be found, for
example, in Chapter~1 of \cite{Olver}.  We note that it follows 
from (\ref{preliminaries:nonoscillatory:5}) that $r_0=1$, and that
the coefficients $r_1,r_2\ldots$ satisfy the recurrence relation
\begin{equation}
r_n = r_{n-1}\left( \frac{\mu- (2n-1)^2 }{4} \right) \frac{2n-1}{2n}.
\label{preliminaries:nonoscillatory:8}
\end{equation}
Using (\ref{preliminaries:nonoscillatory:7}) and (\ref{preliminaries:nonoscillatory:8}),
as many terms as desired in the expansion (\ref{preliminaries:nonoscillatory:alphap_asym})
can be calculated either numerically or analytically.  The first few are
\begin{equation}
\alpha_\nu'(t) \sim 1 - \frac{\mu-1}{8t^2} - \frac{\mu^2-26\mu+25}{128 t^4}
- \frac{\mu^3 -115 \mu^2 +1187\mu - 1073 }{1024t^6} + \cdots.
\label{preliminaries:nonoscillatory:9}
\end{equation}

Obviously, the indefinite integral of  (\ref{preliminaries:nonoscillatory:9})
is
\begin{equation}
\alpha_\nu(t) \sim \tilde{C} + t  + \frac{\mu-1}{8t} + \frac{\mu^2-26\mu+25}{384 t^3}
+ \frac{\mu^3 -115 \mu^2 +1187\mu - 1073 }{5120t^5} + \cdots
\label{preliminaries:nonoscillatory:alpha_asym0}
\end{equation}
with $\tilde{C}$ a constant to be determined to ensure
compatibility with the definition (\ref{preliminaries:nonoscillatory:alpha}).
From the well-known asymptotic expansions
\begin{equation}
J_\nu(t) \sim \sqrt{ \frac{2}{\pi t} } \cos\left(t- \frac{\pi}{2} \nu - \frac{\pi}{4}\right)
\ \ \mbox{as} \ \ t\to\infty
\end{equation}
and
\begin{equation}
Y_\nu(t) \sim \sqrt{ \frac{2}{\pi t} } \sin\left(t- \frac{\pi}{2} \nu - \frac{\pi}{4}\right)
\ \ \mbox{as} \ \ t\to\infty,
\end{equation}
which can be found in Section~7.21 of \cite{Watson} (for example),
we see that   it must be the case that 
\begin{equation}
\tilde{C} = -\frac{\pi}{4} - \frac{\pi}{2}\nu + 2\pi L
\label{preliminaries:nonoscillatory:10}
\end{equation}
with $L$ an arbitrary integer.  In order to set the integer $L$, we evaluated
$\alpha_\nu$ numerically via
(\ref{preliminaries:nonoscillatory:alpha}) for many large values of $t$ and $\nu$
and found that (\ref{preliminaries:nonoscillatory:alpha_asym0}) coincides with $\alpha_\nu$ when
$L$ is taken to be $0$.  That is, 
\begin{equation}
\alpha_\nu(t) \sim -\frac{\pi}{4} - \frac{\pi}{2} \nu + t  + \frac{\mu-1}{8t} + \frac{\mu^2-26\mu+25}{384 t^3}
+ \frac{\mu^3 -115 \mu^2 +1187\mu - 1073 }{5120t^5} + \cdots
\ \ \mbox{as}\ \ t \to \infty.
\label{preliminaries:nonoscillatory:alpha_asym}
\end{equation}

By differentiating (\ref{preliminaries:nonoscillatory:9}) we obtain
\begin{equation}
\alpha_\nu''(t) \sim
   \frac{\mu-1}{4t^3} 
+ \frac{\mu^2-26\mu+25}{32 t^5}
+  \frac{3\mu^3 -345 \mu^2 +3561\mu - 3219 }{512 t^7} + \cdots
\ \ \ \mbox{as}\ \ t \to\infty.
\label{preliminaries:nonoscillatory:alphapp_asym}
\end{equation}
In many instance, differentiating both sides of an asymptotic expansion
such as (\ref{preliminaries:nonoscillatory:9}) will not yield  a valid expression.
However, in this case it is permissible because 
$\alpha_\nu'$ is analytic as function of $t$ 
(see, for instance, Chapter~I of \cite{Olver} for a discussion of this issue).

While the construction of $\alpha_\nu$
presented here is  highly specialized
to the case of Bessel's differential equation,
the existence of nonoscillatory phase functions is an extremely general
phenomenon.  See \cite{Bremer-Rokhlin} for a proof
that under mild conditions on the coefficient $q$ the differential equation
\begin{equation}
y''(t) + q(t) y(t) = 0\ \ \mbox{for all}\ \ a < t < b
\end{equation}
admits one.

\label{preliminaries:nonoscillatory}
\end{subsection}

\begin{subsection}{Univariate Chebyshev series expansions}

For $-1 \leq x \leq 1$ and integers $k \geq 0$, the Chebyshev polynomial $T_k$ 
of degree $k$ is given by the formula 
\begin{equation}
T_k(x) = \cos( k \arccos(x) ).
\label{preliminaries:chebyshev:1}
\end{equation}
The Chebyshev series of a continuous function $f:[-1,1] \to \mathbb{R}$
is 
\begin{equation}
\sideset{}{'}\sum_{j=0}^\infty a_j T_j(x),
\label{preliminaries:chebyshev:2}
\end{equation}
where the coefficients $a_0, a_1, \ldots$ are defined via 
\begin{equation}
a_j = \frac{2}{\pi} \int_{-1}^1 f(x) T_j(x) \frac{dx} { \sqrt{1-x^2} }
\label{preliminaries:chebyshev:3}
\end{equation}
and the dash next to the summation symbol indicates that the first  term in the series
is halved.  It is a consequence
of the well-known relationship between Fourier and Chebyshev
series (as described, for instance, in Chapter~5 of \cite{Mason-Hanscomb})
and the celebrated theorem of Carleson \cite{Carleson}  that 
(\ref{preliminaries:chebyshev:2}) converges to $f$ pointwise almost everywhere on $[-1,1]$.
Similarly, well-known results regarding the convergence of Fourier series
imply that under mild smoothness assumptions on $f$, (\ref{preliminaries:chebyshev:2})
converges uniformly to $f$ on $[-1,1]$.  See Theorem~5.7 in \cite{Mason-Hanscomb},
which asserts that this is the case when $f$ is of bounded variation,
for an example of a result of this type.

If $f:[-1,1] \to\mathbb{R}$ can be analytically continued to the set
\begin{equation}
E_r = \left\{z : \left|z + \sqrt{z^2+1} \right| < r \right\},
\label{preliminaries:chebyshev:5}
\end{equation}
where $r > 1$, then the rate of convergence of (\ref{preliminaries:chebyshev:2})
can be estimated as follows:
\begin{equation}
\sup_{x\in[-1,1]}\left| \sideset{}{'}\sum_{j=0}^N a_j T_j(x) - f(x) \right| = \mathcal{O} \left(r^{-N}\right)
\ \ \mbox{as}\ \ N \to \infty.
\label{preliminaries:chebyshev:6}
\end{equation}
This result can be found as Theorem~5.16 in \cite{Mason-Hanscomb}, among
many other sources.

\label{preliminaries:chebyshev1}
\end{subsection}


\begin{subsection}{Bivariate Chebyshev series expansions}

The bivariate Chebyshev series of a continuous function $f:[-1,1] \times [-1,1] \to \mathbb{R}$
is 
\begin{equation}
\sideset{}{'}\sum_{i=0}^\infty  \sideset{}{'}\sum_{j=0}^\infty a_{ij} T_i(x) T_j(y),
\label{preliminaries:chebyshev:7}
\end{equation}
where the coefficients are defined via the formula
\begin{equation}
a_{ij} = 
\frac{4}{\pi^2}\int_{-1}^1 \int_{-1}^1 f(x,y) T_i(x) T_j(y) \frac{dx}{\sqrt{1-x^2}}\frac{dy}{\sqrt{1-y^2}}
\label{preliminaries:chebyshev:8}
\end{equation}
and the dashes next to the summation symbols indicate that the first  term in each
sum is halved.    It is an immediate consequence  of the results of 
\cite{Fefferman1} on the pointwise almost everywhere convergence of multiple Fourier series 
 that 
\begin{equation}
\lim_{N\to\infty} \sideset{}{'}\sum_{i=0}^N \sideset{}{'}\sum_{j=0}^N
a_{ij} T_i(x) T_j(y)  = f(x,y)
\label{preliminaries:chebyshev:9}
\end{equation}
for almost all $(x,y) \in [-1,1] \times [-1,1]$.   
As in the case of univariate Chebyshev series,
under mild smoothness conditions on $f$, the convergence of (\ref{preliminaries:chebyshev:9}) 
is uniform.  See, for instance, Theorem~5.9 in  \cite{Mason-Hanscomb}.

The result of the preceding section on the convergence of the Chebyshev series
of analytic functions can be generalized to bivariate Chebyshev expansions.
In particular, if  $f(x,y)$ is analytic on the set 
\begin{equation}
\left\{
\left(x,y\right) \in \mathbb{C} \times \mathbb{C} :\
\left|x+\sqrt{x^2-1}\right| < r_1,\  \ 
\left|y+\sqrt{y^2-1}\right| < r_2 
\right\},
\label{preliminaries:chebyshev:10}
\end{equation}
where $r_1, r_2 >1$, then
\begin{equation}
\sup_{x\in[-1,1]}\left|\sideset{}{'}\sum_{i=0}^N \sideset{}{'}\sum_{j=0}^N
a_{ij} T_i(x) T_j(x) - f(x,y) \right|
= O \left(r_1^{-N} r_2^{-N} \right) \ \ \mbox{as}\ \ N \to \infty.
\label{preliminaries:chebyshev:11}
\end{equation}
This follows from Theorem~11 in Chapter~V of \cite{Bochner-Martin}.

\label{preliminaries:chebyshev2}
\end{subsection}


\begin{subsection}{Chebyshev interpolation and spectral integration}

In practice it is, of course, not possible to compute all of the coefficients
in the Chebyshev series expansions (\ref{preliminaries:chebyshev:2}) or
(\ref{preliminaries:chebyshev:7}), or even to compute the first few coefficients
in these series exactly.  Standard results, however,
show that these coefficients of such a series can be approximated to high 
accuracy assuming that they decay rapidly.

For each nonnegative integer $n$, we refer to the collection of points
\begin{equation}
\rho_{j,n} =  -\cos\left(\frac{\pi j}{n}\right), \ \ j=0,1,\ldots,n,
\label{preliminaries:chebyshev:nodes}
\end{equation}
as the $(n+1)$-point Chebyshev grid on the interval $[-1,1]$,
and we call individual elements of this set  Chebyshev nodes or points.
One discrete version of  the well-known orthogonality relation
\begin{equation}
\int_{-1}^1   \frac{T_i(x) T_j(x)}{\sqrt{1-x^2}}\ dx = 
\begin{cases}
0   & \mbox{if}\ \ i \neq j \\
\frac{\pi}{2} & \mbox{if}\ \ i = j > 0 \\
\pi & \mbox{if}\ \ i=j=0.
\end{cases}
\end{equation}
is
\begin{dmath}
\sideset{}{''}\sum_{l=0}^{n} T_i(\rho_{l,n}) T_j(\rho_{l,n})
=
\begin{cases}
0   & \mbox{if}\ \ 0 \leq i,j \leq  n\ \ \mbox{and}\ i \neq j \\
\frac{n}{2} & \mbox{if}\ \ 0 < i=j <n\\
n & \mbox{if}\ \ i=j=0\ \  \mbox{or}\ \ i=j=n.
\end{cases}
\label{preliminaries:chebyshev:12}
\end{dmath}
Here, the double dash next to the summation sign
indicates that the first and last term in the series are halved.
Formula~(\ref{preliminaries:chebyshev:12})
 can be found in a slightly different form in Chapter~4 of \cite{Mason-Hanscomb}.
Any univariate polynomial $f$ of degree $n$ can be represented
in the form
\begin{equation}
f(x) = \sideset{}{''}\sum_{l=0}^n b_l T_l(x),
\label{preliminaries:chebyshev:13}
\end{equation}
and the coefficients in the expansion (\ref{preliminaries:chebyshev:13})
can be easily computed  from the values of $f$ at the nodes 
(\ref{preliminaries:chebyshev:nodes}).  In particular,  it follows
from (\ref{preliminaries:chebyshev:12})  that
\begin{equation}
b_k = \frac{2}{n} \sideset{}{''}\sum_{l=0}^n T_k\left(\rho_{l,n}\right) f\left(\rho_{l,n}\right)
\label{preliminaries:chebyshev:14}
\end{equation}
for all $k=0,1,\ldots,n$.

If $f:[-1,1] \to \mathbb{R}$ is smooth but no longer a polynomial, then the coefficients 
$b_0,b_1,\ldots,b_n$ obtained from (\ref{preliminaries:chebyshev:14})
are related to the coefficients $a_0,a_1,\ldots$ in the Chebyshev expansion
(\ref{preliminaries:chebyshev:2}) of $f$ via 
\begin{equation}
b_k = a_k + \sum_{j=1}^\infty  \left( a_{k+2jn} + a_{-k+2jn} \right)
\ \ \mbox{for all} \ \ k=0,1,\ldots,n.
\label{preliminaries:chebyshev:15}
\end{equation}
 This result can be found in a slightly
different form in Section~6.3.1 of \cite{Mason-Hanscomb}.  It follows
easily from (\ref{preliminaries:chebyshev:15}) and the fact that
the Chebyshev polynomials are bounded in $L^\infty\left([-1,1]\right)$
norm by $1$ that
\begin{equation}
\sup_{x \in [-1,1]} 
\left| f(x) - \sideset{}{''}\sum_{l=0}^n b_l T_l(x)  \right|
\leq 
2 \sum_{l={n+1}}^\infty |a_l|.
\label{preliminaries:chebyshev:16}
\end{equation}
In other words, assuming that the coefficients of the Chebyshev expansion
of $f$ decay rapidly, the series (\ref{preliminaries:chebyshev:13})
converges rapidly to $f$ as $n \to \infty$.  
We will, by a slight abuse of terminology,  refer to (\ref{preliminaries:chebyshev:13}) 
as the $n^{th}$ order Chebyshev 
expansion of the function $f:[-1,1] \to \mathbb{R}$.

Given only the values of a function $f:[-1,1] \to \mathbb{R}$
at the nodes of the $(n+1)$-point Chebyshev grid on $[-1,1]$, it is possible
to evaluate the Chebyshev  expansion (\ref{preliminaries:chebyshev:13})
in an efficient and  numerical stable fashion  without explicitly
computing the coefficients $b_0,\ldots,b_n$.   In particular, the 
value of  (\ref{preliminaries:chebyshev:13}) at a  point $x \in [-1,1]$
which does not coincide with any of the grid points $\rho_{0,n},\ldots,\rho_{n,n}$
is given by the barycentric interpolation formula
\begin{equation}
\left(\sum_{j=0}^n \frac{(-1)^j f(\rho_{j,n})}{x-\rho_{j,n}}\right)
\Big/
\left(\sum_{j=0}^n \frac{(-1)^j}{x-\rho_{j,n}}\right).
\label{preliminaries:chebyshev:barycentric}
\end{equation}
See, for instance, \cite{Trefethen} for a thorough discussion of the numerical
stability and efficiency of this technique.

For each $k > 1$, the general antiderivative of $T_k$ is 
\begin{equation}
\int T_k(x)\ dx = 
\frac{1}{2} \left( \frac{T_{k+1}(x)}{k+1} - \frac{T_{k-1}(x)}{k-1} \right) + C
\end{equation}
while
\begin{equation}
\int T_0(x)\ dx = \frac{1}{4} T_1(x) + C
\end{equation}
and
\begin{equation}
\int T_1(x)\ dx = \frac{1}{4} T_2(x) + C.
\end{equation}
These formulas can be found, for instance,
 in Section~2.4.4 of \cite{Mason-Hanscomb}.  Using them, the values
of 
\begin{equation}
g(x) = \int_{-1}^x  f(t)\ dt,
\end{equation}
where  $f$ is defined via (\ref{preliminaries:chebyshev:13}), can be 
computed at the Chebyshev nodes  (\ref{preliminaries:chebyshev:nodes}).
We will refer to the matrix which takes the values
of $f$ at the nodes (\ref{preliminaries:chebyshev:13}) to those
of $g$ at the same nodes as the $(n+1) \times (n+1)$ spectral
integration matrix.

There results for univariate Chebyshev expansions can be easily generalized
to the case of bivariate Chebyshev expansions.  Given
$f:[-1,1]\times [-1,1]\to\mathbb{R}$, we will refer to the
series
\begin{equation}
\sideset{}{''} \sum_{i=0}^n \sideset{}{''} \sum_{j=0}^n b_{ij} T_i(x) T_j(y)
\label{preliminaries:chebyshev:17}
\end{equation}
whose  the coefficients $\{b_{ij} : 0 \leq i, j \leq n\}$ are defined via
the formula
\begin{equation}
b_{ij} = 
\frac{4}{n^2}
\sideset{}{''} \sum_{l=0}^n  \sideset{}{''}\sum_{k=0}^n 
T_i\left(\rho_{l,n}\right) T_j\left(\rho_{k,n}\right) f\left(\rho_{l,n}, \rho_{k,n}\right)
\label{preliminaries:chebyshev:18}
\end{equation}
as the $n^{th}$ order Chebyshev expansion of $f$.  It is easy to see that
\begin{equation}
\sup_{x \in [-1,1]}
\left|
f(x,y) - \sideset{}{''} \sum_{i=0}^n \sideset{}{''} \sum_{j=0}^n b_{ij} T_i(x) T_j(y)
\right|
\leq  2 \sum_{i=n+1}^\infty \sum_{j=n+1}^\infty \left|a_{ij}\right|,
\end{equation}
where the $\{ a_{ij}\}$ are defined via (\ref{preliminaries:chebyshev:8}).

\label{preliminaries:chebyshev3}
\end{subsection}


\begin{subsection}{Compressed bivariate Chebyshev expansions}

It often happens that  many of the coefficients
in the bivariate Chebyshev expansion
(\ref{preliminaries:chebyshev:17}) of a function $f:[-1,1]\times [1,1] \to \mathbb{R}$
are of negligible magnitude.  In order to reduce the cost of storing 
 such expansions as well as the cost of evaluating them, we use the following
construction to reduce the number of coefficients which need to be considered.

Suppose that $\epsilon > 0$, and that 
\begin{equation}
\sideset{}{''}\sum_{i=0}^n \sideset{}{''}\sum_{j=0}^{n} b_{ij} T_i(x) T_j(y)
\label{preliminaries:compressed:1}
\end{equation}
is the $n^{th}$ order Chebyshev expansion for $f:[-1,1] \times [-1,1] \to \mathbb{R}$.
We let $M$ denote the least positive integer which is less than
or equal to $n$ and such that
\begin{equation}
\left|b_{ij}\right| < \epsilon
\ \ \mbox{for all} \ \ i > M\ \mbox{and} \ \   j=0,\ldots,n,
\end{equation}
assuming such an integer exists.  If it does not, then we take $M=n$.
Similarly, for each $i=0,\ldots,M$, we let $m_i$ be the least positive integer
less than or equal to $n$ such that
\begin{equation}
\left|b_{ij}\right| < \epsilon
\ \ \mbox{for all}\ \  j=m_i+1,\ldots,n
\end{equation}
if such an integer exists, and we let $m_i = n$ otherwise.
We refer to the series
\begin{equation}
\sum_{i=0}^M\sum_{j=0}^{m_i} \widetilde{b_{ij}} T_i(x) T_j(x),
\label{preliminaries:chebyshev4:expansion}
\end{equation}
where $\widetilde{b_{ij}}$ is defined via
\begin{equation}
\widetilde{b_{ij}} =
\begin{cases}
b_{ij} & \mbox{if} \  \ 1 \leq i,j \leq n \\
\frac{1}{2} b_{ij} & \mbox{if} \  \ 1 \leq j \leq n \ \ \mbox{and}\ \ i = 0,\ n\\
\frac{1}{2} b_{ij} & \mbox{if} \  \ 1 \leq i \leq n \ \ \mbox{and}\ \ j = 0,\ n\\
\frac{1}{4} b_{ij} & \mbox{otherwise},
\end{cases}
\end{equation}
as the $\epsilon$-compressed $n^{th}$ order Chebyshev expansion of $f$.

Obviously, the results discussed in Sections~\ref{preliminaries:chebyshev1} 
through Section~\ref{preliminaries:chebyshev4} can be modified in 
a straightforward fashion so as to apply to a function given on an arbitrary 
interval $[a,b]$ (in the case of univariate functions) 
or one given on a compact rectangle $[a,b] \times [c,d]$ (in the case  
of a bivariate function).   For instance,   the nodes of the $(n+1)$-point Chebyshev grid on $[a,b]$
are
\begin{equation}
\tilde{\rho}_{j,n} = - \frac{b-a}{2}\cos\left(\frac{\pi j}{n} \right) + \frac{b+a}{2},
\ \ \ j=0,1,\ldots,n
\end{equation}
and the $n^{th}$ order Chebyshev expansion of $f:[a,b] \to \mathbb{R}$ is 
\begin{equation}
\sideset{}{''}\sum_{i=0}^n b_i T_{i} \left( \frac{2}{b-a} x + \frac{b+a}{a-b} \right),
\end{equation}
where the coefficients $b_0,\ldots,b_n$ are given by
\begin{equation}
b_i = \frac{2}{n} \sideset{}{''}\sum_{l=0}^n T_i\left(\rho_{l,n}\right) 
f\left(\tilde{\rho}_{l,n}\right).
\end{equation}

\label{preliminaries:chebyshev4}
\end{subsection}


\begin{subsection}{An adaptive discretization procedure}

We now briefly describe a fairly standard procedure for
 adaptively discretizing a smooth function $f:[a,b] \to \mathbb{R}$.
It takes as input a desired precision $\epsilon >0$ and a positive integer $n$.
The goal of this procedure is to construct a partition
\begin{equation}
a = \gamma_0 < \gamma_1 < \cdots < \gamma_m = b
\end{equation}
of $[a,b]$ such that  the $n^{th}$ order Chebyshev expansion of $f$ on each of the subintervals
$[\gamma_j,\gamma_{j+1}]$ of $[a,b]$ approximates $f$ with accuracy $\epsilon$.  That is,
for each $j=0,\ldots,m-1$ we aim to achieve
\begin{equation}
\sup_{x \in [\gamma_j,\gamma_{j+1}]}
\left| 
f(x) - 
\sideset{}{''}\sum_{i=0}^n b_{i,j} T_{i} \left( \frac{2}{\gamma_{j+1}-\gamma_j } x + 
\frac{\gamma_{j+1}+\gamma_j}{\gamma_j-\gamma_{j+1}} \right)
\right| < \epsilon,
\label{preliminaries:adaptive:1}
\end{equation}
where $b_{0,j},b_{1,j}\ldots,b_{n,j}$ are the coefficients in the $n^{th}$ order Chebyshev expansion
of $f$ on the interval $\left[\gamma_j,\gamma_{j+1}\right]$.  
These coefficients are defined by the formula
\begin{equation}
b_{i,j} = \frac{2}{n} \sideset{}{''}\sum_{l=0}^n T_i\left(\rho_{l,n}\right) 
f\left( \frac{\gamma_{j}-\gamma_{j+1}}{2}\cos\left(\frac{\pi l}{n} \right) + \frac{\gamma_{j+1}+\gamma_j}{2}\right).
\end{equation}

During the procedure,  two lists of subintervals are maintained: a list
of subintervals which are to be processed and a list of output subintervals.
Initially, the list of subintervals to be processed consists of $[a,b]$ 
and the list of output subintervals is empty.
The procedure terminates when the list of subintervals to be processed is empty
or when the number of subintervals in this list  exceeds a present
limit (we usually take this limit to be $300$).
In the latter case, the procedure is deemed to have failed.
As long as the list of subintervals to process is nonempty and its length
does not exceed the preset maximum, the algorithm
proceeds by removing a subinterval $\left[\eta_1,\eta_2\right]$
from that list  and performing the following operations:
\begin{enumerate}

\item
Compute the coefficients $b_0,\ldots,b_n$ 
in the $n^{th}$ order Chebyshev expansion
of the restriction of $f$ to the interval $\left[\eta_1,\eta_2\right]$.

\vskip 1em
\item
Compute the quantity 
\begin{equation}
\Delta = 
\frac{\max\left\{
\left|b_{\frac{n}{2}+1} \right|, 
\left|b_{\frac{n}{2}+2} \right|, \ldots
\left|b_n \right|
\right\}}
{\max\left\{
\left|b_0 \right|, 
\left|b_1 \right|, \ldots
\left|b_n \right|
\right\}}.
\end{equation}

\vskip 1em
\item
If $\Delta < \epsilon$ then 
 the subinterval $\left[\eta_1,\eta_2\right]$ is added to the list
of output subintervals.

\vskip 1em
\item 

If $\Delta \geq \epsilon$, then the subintervals
\begin{equation}
\left[\eta_1,\frac{\eta_1+\eta_2}{2}\right]
\ \ \mbox{and}\ \ 
\left[\frac{\eta_1+\eta_2}{2}, \eta_2 \right]
\end{equation}
are added to the list of  subintervals to be processed.

\end{enumerate}

This algorithm is  heuristic in the sense that there is no guarantee
that (\ref{preliminaries:adaptive:1}) will be achieved, but 
similar adaptive discretization procedures are widely used
with great success.  However, there is one common circumstance which 
leads to the failure of this procedure.  The quantity $\Delta$ is an attempt to estimate
the relative accuracy with which the Chebyshev expansion
of $f$ on the interval $\left[\eta_1,\eta_2\right]$ approximates
$f$.  In cases in which the condition number of the evaluation of
$f$ is larger than $\epsilon$ on some part of $[a,b]$,
the procedure will generally fail or an excessive number of subintervals
will be generated.   Particular care needs to be taken when $f$ has a zero
in $[a,b]$.
In most cases, for $x$ near  a zero of $f$, the condition number of evaluation
of $f(x)$ (as defined in Section~\ref{preliminaries:condition})
is large.   In this article, we avoid such difficulties  by always
considering functions which are bounded away from $0$.

\label{preliminaries:adaptive}
\end{subsection}

\label{section:preliminaries}
\end{section}

\begin{section}{An adaptive solver for nonlinear differential equations}

In this section, we describe a numerical algorithm for the solution
of nonlinear second order differential equations of the form
\begin{equation}
y''(t) = f(t,y(t),y'(t)) \ \ \mbox{for all}\ \  a < t < b.
\label{solver:1}
\end{equation}
It is intended to be extremely robust, but not necessarily
highly efficient.  This is fitting since we only use it to perform precomputations.
Here, we assume that initial conditions for the desired solution
$y$ of (\ref{solver:1}) are specified.   That is, we seek a solution of (\ref{solver:1})
which satisfies
\begin{equation}
y(a) = y_a \ \ \mbox{and} \ \ y'(a) = y_a',
\label{solver:2}
\end{equation}
where the constants $y_a$ and $y_a'$ are given.  The algorithm can easily be 
modified for the case of a terminal value  problem.

The procedure takes as input a subroutine 
for evaluating $f$ and its derivatives with respect to $y$ and $y'$,
a positive integer
$n$, and a positive real number $\epsilon >0$.  It maintains
a list of output subintervals, a stack 
containing a set of subintervals to process, and two constants $c_1$ and $c_2$.
Initially, the list of output subintervals is empty,
the stack  consists of $[a,b]$, $c_1$ is taken
to be $y_a$, and $c_2$ is taken to be $y_a'$. 
If the size of the stack exceeds a preset maximum (taken to be
$300$ in the calculations performed in this paper), the procedure is deemed
to have failed.  As long as the stack is nonempty and its length
does not exceed the preset maximum, the algorithm
proceeds by popping an interval $\left[\eta_1,\eta_2\right]$
off of the stack and performing the following operations:
\begin{enumerate}

\item
Form the nodes $t_0, t_1,\ldots,t_n$ of the $(n+1)$-point Chebyshev grid
on the interval $\left[\eta_1,\eta_2\right]$.

\vskip 1em
\item
Apply the trapezoidal method (see, for instance, \cite{Iserles}) in order to approximate the 
values of the second derivative $y_0''$ of the function satisfying 
\begin{equation}
\left\{
\begin{aligned}
y_0''(t) &= f(t,y(t),y'(t)) \ \  \mbox{for all} \ \ \eta_1 \leq t \leq \eta_2 \\
y_0(\eta_1) &= c_1 \\
y_0'(\eta_1) &= c_2
\end{aligned}
\right.
\label{solver:3}
\end{equation}
at the points $t_0, t_1, \ldots, t_{n}$.  

\vskip 1em
\item

Using spectral integration,  approximate the values of $y_0'$ and $y_0$ 
at the points $t_0,\ldots,t_n$ through the formulas
\begin{equation}
y_0'(t) = c_2 + \int_{\eta_1}^t y_0''(s)\ ds
\label{solver:4}
\end{equation}
and
\begin{equation}
y_0(t) = c_1 + \int_{\eta_1}^t y_0'(s)\ ds.
\label{solver:5}
\end{equation}

\vskip 1em
\item
Apply Newton's method to the initial value problem (\ref{solver:3}).
The function $y_0$ is used as the initial guess for Newton's method.
The $j^{th}$ iteration of Newton's method starts with an initial
approximation $y_{j-1}$ and consists of solving the linearized problem 
\begin{equation}
 \delta_j''(t)
-\frac{\partial f}{\partial y}(t,y_0(t),y_0'(t)) \delta_j(t) 
-
\frac{\partial f}{\partial y'}(t,y_0(t),y_0'(t)) \delta_j'(t) 
 = f(t,y_0(t),y_0'(t)) - y_0''(t)
\label{solver:6}
\end{equation}
on the interval $\left[\eta_1,\eta_1\right]$
for $\delta_j$ and forming the new approximation  $y_j = y_{j-1} + \delta_j$.
The initial conditions 
\begin{equation}
\delta_j(0) = \delta_j'(0) = 0
\label{solver:7}
\end{equation}
are imposed 
since the initial approximation $y_0$ is already consistent with
the desired initial conditions.
All of the functions appearing in (\ref{solver:6}) are represented 
via their values at the points $t_0,t_1,\ldots,t_n$.

\vskip 1em
An integral equation method is used to solve (\ref{solver:6}).
More specifically, by assuming  that $\delta_j$ is given by 
\begin{equation}
\delta_j(t) =
\int_{\eta_1}^t \int_{\eta_1}^s \sigma_j(\tau)\ d\tau \ ds 
\label{solver:8}
\end{equation}
the system (\ref{solver:6})  is transformed into the integral equation
\begin{equation}
\begin{aligned}
\sigma_j(t) 
-&\frac{\partial f}{\partial y}(t,y_{j-1}(t),y_{j-1}'(t)) 
\int_{\eta_1}^t \int_{\eta_1}^s \sigma_j(\tau)\ d\tau \ ds 
-\\
&\frac{\partial f}{\partial y'}(t,y_{j-1}(t),y_{j-1}'(t)) 
 \int_{\eta_1}^t \sigma_j(\tau)\ d\tau \ ds 
  = f(t,y_{j-1}(t),y_{j-1}'(t)) - y_{j-1}''(t).
\end{aligned}
\label{solver:10}
\end{equation}
Note that the choice of the  representation 
(\ref{solver:8}) of $\delta_j$ is consistent with the  conditions (\ref{solver:7}).
The linear system which arises from requiring that
(\ref{solver:10}) is satisfied at the points $t_0,t_1,\ldots,t_n$
is inverted in order to calculate the values of $\sigma_j$
at those nodes.
Spectral integration (as described in Section~\ref{preliminaries:chebyshev4}) is used to evaluate
the integrals, and to compute the values of $\delta$ at $t_0,\ldots,t_n$ 
via (\ref{solver:8}).
\vskip 1em

Define
\begin{equation}
\Delta_j = \max\left\{\left|\delta_j(t_0)\right|,
\left|\delta_j(t_1)\right|,\ldots,\left|\delta_j(t_n)\right|\right\}.
\end{equation}
If $j > 1$ and   $\Delta_j < \Delta_{j-1}$, then Newton iterations continue.
Otherwise, the Newton iteration is terminated, having obtained $y_{j-1}$ 
as the result of the procedure.

\vskip 1em
\item
Compute the Chebyshev coefficients $b_0,b_1,\ldots,b_n$ 
of the polynomial which interpolate the values
of $y_{j-1}$ at the points $t_0,t_1,\ldots,t_n$
and  define $\Lambda$ via
\begin{equation}
\Lambda = 
\frac{\max\left\{\left|b_{\frac{n}{2}+1}\right|,\left|b_{\frac{n}{2}+2}\right|, \ldots \left|b_{n}\right|\right\}}
{\max\left\{\left|b_1\right|,\left|b_{2}\right|, \ldots \left|b_n\right|\right\}}.
\label{solver:11}
\end{equation}
If $\Lambda > \epsilon$, then push the subintervals
\begin{equation}
\left[\frac{\eta_1+\eta_2}{2},\eta_2\right]
\end{equation}
and 
\begin{equation}
\left[\eta_1,\frac{\eta_1+\eta_2}{2}\right]
\label{solver:1000}
\end{equation}
 onto the stack (in that order) so that  (\ref{solver:1000}) is the next interval
to be processed by the algorithm.

\vskip 1em
If $\Lambda \leq \epsilon$, then 
 $\left[\eta_1,\eta_2\right]$ is added to the list of output intervals,
$c_1$ is set equal to $y_j(\eta_2)$,  and $c_2$ is set equal to $y_j'(\eta_2)$.
\end{enumerate}

\vskip 1em

As in the case of the adaptive discretization procedure of 
Section~\ref{preliminaries:adaptive}, this is a heuristic
algorithm which is not guaranteed to achieve an accurate
discretization of the solution of (\ref{solver:1}).
Moreover, the quantity $\Lambda$ defined in  (\ref{solver:11}) is an attempt
to measure the relative accuracy with which the obtained solution 
of (\ref{solver:1}) is represented on the interval under consideration.
When the condition number of evaluation of the solution $y$ of (\ref{solver:1}) 
is large, this algorithm tends to produce an excessive number
of intervals or fail altogether.  Since the condition number of evaluation
 of a function $f$ is generally large near its zeros, in this article
we always apply it in  cases in which the solution of (\ref{solver:1})
is bounded away from $0$.

\vskip 1em
\begin{remark}
When applied to (\ref{solver:1}), the trapezoidal method
produces approximations of the values of $y_0$ and $y_0'$ 
at the nodes $t_0,t_1,\ldots,t_n$ in addition to approximations
of the values of $y_0''$.  We discard those values and recompute
$y_0$ and $y_0'$ via spectral integration.    We do so
because while the values of $y_0$ and $y_0'$ obtained
by the trapezoidal method at $t_0,t_1,\ldots,t_n$ must
satisfy the realtions
\begin{equation*}
y_0''(t_j) = f(t_j,y(t_j),y'(t_j)) \ \  \mbox{for all} \ \ j=0,1,\ldots,n,
\end{equation*}
they need not be consistent with each other in the sense that
\begin{equation}
\left(
\begin{array}{c}
y_0'(t_0) \\
y_0'(t_1) \\
\vdots\\
y_0'(t_n)
\end{array}
\right)
=
\left(
\begin{array}{c}
c_2 \\
c_2 \\
\vdots\\
c_2
\end{array}
\right)
+
S_n
\left(
\begin{array}{c}
y_0''(t_0) \\
y_0''(t_1) \\
\vdots\\
y_0''(t_n)
\end{array}
\right)
\end{equation}
and
\begin{equation}
\left(
\begin{array}{c}
y_0(t_0) \\
y_0(t_1) \\
\vdots\\
y_0(t_n)
\end{array}
\right)
=
\left(
\begin{array}{c}
c_1 \\
c_1 \\
\vdots\\
c_1
\end{array}
\right)
+
S_n
\left(
\begin{array}{c}
y_0'(t_0) \\
y_0'(t_1) \\
\vdots\\
y_0'(t_n)
\end{array}
\right),
\end{equation}
where $S_n$ denotes the $(n+1)\times (n+1)$ spectral integration matrix
might not hold.  Proceeding without recomputing the values of $y_0$
and $y_0'$ in order to make sure that these consistency conditions
are satisfied would lead to the failure of  Newton's method in most
cases.
\end{remark}

\label{section:solver}
\end{section}

\begin{section}{An algorithm for the rapid numerical solution of Bessel's differential equation}

In this section, we describe a numerical algorithm for the solution of Bessel's
differential equation  for a fixed value of  $\nu$.
Our algorithm runs in time independent of $\nu$ and is a key component
of the scheme of the following section for the construction of 
tables which allow for the rapid numerical evaluation of the Bessel functions.

The algorithm takes as input $\nu >0$ and a desired precision $\epsilon > 0$.
It proceeds in three stages.

{\it Stage one: computation of a nonoscillatory phase function}

In this stage, we calculate
the nonoscillatory phase function $\alpha_\nu$ defined by (\ref{preliminaries:nonoscillatory:alpha})
on the interval 
\begin{equation}
\left[ \sqrt{\nu^2-\frac{1}{4}},\ 1000\ \nu\right]
\end{equation}
if $\nu > \frac{1}{2}$, and on the interval
\begin{equation}
\left[ 2, 1000 \right]
\end{equation}
in the event that  $0 \leq \nu \leq \frac{1}{2}$.  In either case, we will denote 
the left-hand side of the interval on which we calculate $\alpha_\nu$  by $a$ and 
the right-hand side by $b$.  

We first construct the derivative $\alpha_\nu'$ of $\alpha_\nu$ by solving 
 Kummer's equation (\ref{preliminaries:riccati:kummer}) on the interval
$[a,b]$ with  with $q$ taken to be the coefficient of $y$ in Bessel's differential
equation (\ref{introduction:besseleq}); that is,
\begin{equation}
q(t) = 1 - \frac{\nu^2 - \frac{1}{4}}{t^2}.
\label{phase:q}
\end{equation}
Most solutions of (\ref{preliminaries:riccati:kummer}) are oscillatory;
however, the phase function $\alpha_\nu$ is a nonoscillatory.
 Moreover,  the values of $\alpha_\nu'$ and
its derivative $\alpha_\nu''$ at the right-hand endpoint $b$
can be approximated to high accuracy via the asymptotic expansions 
(\ref{preliminaries:nonoscillatory:alphap_asym}) and 
(\ref{preliminaries:nonoscillatory:alphapp_asym}).  
Accordingly, we solve a terminal value problem for Kummer's equation
with the values of $\alpha_\nu'(b)$ and $\alpha_\nu''(b)$ specified.
We   use the  adaptive procedure of Section~\ref{section:solver}
to solve Kummer's equation;
the input $n$ to that procedure is taken to be $30$ and 
the desired precision is set to $\epsilon$.
The functions $\alpha_\nu'$ and $\alpha_\nu''$  are represented
via their values at the nodes of the $31$-point Chebyshev grids 
(see Section~\ref{preliminaries:chebyshev3}) 
on a collection of subintervals
\begin{equation}
\left[\gamma_0,\gamma_1\right],\left[\gamma_1,\gamma_2\right],\ldots,\left[\gamma_{m-1},\gamma_m\right],
\label{phase:subintervals}
\end{equation}
where $a = \gamma_0 < \gamma_1 < \ldots < \gamma_m=b$ is a partition of $[a,b]$
which is  determined adaptively by the solver of Section~\ref{section:solver}.
We use the formula
\begin{equation}
\alpha_\nu(t) = \alpha_\nu(b) + \int_b^t \alpha_\nu'(s)\ ds
\label{phase:int}
\end{equation}
to calculate $\alpha_\nu$.   More specifically, spectral integration is used to obtain
the values of $\alpha_\nu$ at the nodes of the $31$-point Chebyshev grids 
on the subintervals (\ref{phase:subintervals}).
The value of $\alpha_\nu(b)$  is approximated
to high accuracy via the asymptotic expansion (\ref{preliminaries:nonoscillatory:alpha_asym}).
We use the first $30$ terms of each of the expansions
(\ref{preliminaries:nonoscillatory:alpha_asym}),
(\ref{preliminaries:nonoscillatory:alphap_asym})
and
(\ref{preliminaries:nonoscillatory:alphapp_asym}).

The functions $\alpha_\nu$, $\alpha_\nu'$ and $\alpha_\nu''$  
can  calculated in an efficient and numerically stable
fashion at any point in the interval $[a,b]$ via
barycentric Chebyshev interpolation using their values
at the nodes of the Chebyshev grids on the subintervals
(\ref{phase:subintervals}) (see Section~\ref{preliminaries:chebyshev3}).
Using the values of $\alpha_\nu$ and $\alpha_\nu'$, 
$J_\nu$ and $Y_\nu$ can be evaluated at any point on the interval
$[a,b]$ via (\ref{preliminaries:nonoscillatory:1})
and (\ref{preliminaries:nonoscillatory:2}).

\vskip 1em
{\it  Step two: computation of $\nu + \log\left(-Y_\nu(t) \sqrt{t}\right)$}

In the event that $\nu > \frac{1}{2}$, we calculate
 the function $\nu + \log\left(-Y_\nu(t)\sqrt{t}\right)$ on the interval
\begin{equation}
\left[ \frac{\nu}{1000}, \sqrt{\nu^2-\frac{1}{2}} \right]
\end{equation}
by solving a terminal value problem for Riccati's equation 
\begin{equation}
r''(t) + (r'(t))^2 + q(t) = 0
\label{phase:riccati}
\end{equation}
with $q$ given by (\ref{phase:q}).  In fact,
we solve Riccati's equation on the slightly larger interval
\begin{equation}
\left[ \frac{\nu}{1000}, t^* \right],
\label{phase:50}
\end{equation}
where $t^*$ is the solution of the equation
\begin{equation}
\alpha_\nu\left(t^*\right) =    -\frac{\pi}{2}.
\label{phase:100}
\end{equation}
That there exists a solution $t^*$ of this equation such that
\begin{equation}
t^* > \sqrt{\nu^2-\frac{1}{4}}
\end{equation}
 is a consequence of a well-known result regarding the zeros of Bessel functions;
namely, that  $J_\nu$ cannot have zeros on the interval
\begin{equation}
\left(0,\sqrt{\nu^2-\frac{1}{4}}\right]
\label{phase:102}
\end{equation}
(see, for instance, Chapter~15 of \cite{Watson}).
From (\ref{preliminaries:nonoscillatory:3}), we see that the zeros of $J_\nu$
occur at points $t$ such that
\begin{equation}
\alpha_\nu(t) = \frac{\pi}{2} + \pi k \ \ \mbox{with}  \ \ k \in \mathbb{Z}.
\label{phase:103}
\end{equation}
It is obvious from the  definition (\ref{preliminaries:nonoscillatory:alpha})
of $\alpha_\nu$ that  $\alpha_\nu(0) = -\frac{\pi}{2}$, and that $\alpha_\nu$ is increasing
as a function of $t$.  Consequently, if $t^*$ denotes the smallest
positive real number such that $J_\nu\left(t^*\right) = 0$, then 
$t^*$ satisfies (\ref{phase:100}) and  and it must be the case that
\begin{equation}
t^* > \sqrt{\nu^2-\frac{1}{4}}
\label{phase:104}
\end{equation}
since there are no zeros of $J_\nu$  in (\ref{phase:102}).
The values of $\alpha_\nu$ and its derivative having been calculated
in the preceding phase, there is no difficulty in using
Newton's method to obtain the value of $t^*$ by solving
the nonlinear equation (\ref{phase:100}) numerically.
Moreover, the values of the functions $Y_\nu$ and $Y_\nu'$ 
at $t^*$ can be calculated without the loss of precision 
indicated by their  condition numbers of evaluation 
(see Section~\ref{preliminaries:condition} for a definition
of the condition number of evaluation of a function).
In particular, since $\alpha_\nu(t^*) = -\frac{\pi}{2}$, 
\begin{equation}
Y_\nu(t^*) = \frac{\sin(\alpha_\nu(t^*))}{\alpha_\nu'(t^*)}
= 
- \frac{1}{\alpha_\nu'(t^*)}
\label{phase:105}
\end{equation}
and
\begin{equation}
Y_\nu'(t^*) 
= \cos(\alpha_\nu(t^*)) \sqrt{\alpha_\nu'(t^*)} 
 - 
\frac{\sin(\alpha_\nu(t^*)) \alpha_\nu''(t^*)}{2 (\alpha_\nu'(t^*))^{\frac{3}{2}} }
=
\frac{\alpha_\nu''(t)}{2 (\alpha_\nu'(t))^{\frac{3}{2}}}.
\label{phase:106}
\end{equation}
The condition number of the evaluation of the nonoscillatory
 functions $\alpha_\nu'$ and $\alpha_\nu''$ is not large and is bounded
independent of  $\nu$,
so there calculations can be performed without much loss of accuracy.
See, for instance, \cite{BremerZeros}, where this issue is  discussed in detail.  
We note that the numerical evaluation of $J_\nu$ and $Y_\nu$ at an arbitrary point
$t$ via (\ref{preliminaries:nonoscillatory:3}) and 
(\ref{preliminaries:nonoscillatory:4}) will result in a relative error
on the order of the condition number of the evaluation of these functions.
This loss of accuracy stems from the evaluation of the trigonometric
functions cosine and sine which appear in those formulas.

From the values of $Y_\nu$ and $Y_\nu'$ at $t^*$, we calculate
the values of $\nu + \log(-Y_\nu(t)\sqrt{t})$ and its derivative there.
Then we solve the corresponding terminal value problem for Riccati's equation.
We use the solver described in Section~\ref{section:solver} to do so.
Our motivation for calculating
$\nu + \log(-Y_\nu(t)\sqrt{t})$
in lieu of $\log(-Y_\nu(t)\sqrt{t})$ is that the former is bounded away from 
$0$ on the interval (\ref{phase:50}) while the latter is not.  As discussed
in Section~\ref{preliminaries:condition} the condition number of evaluation
of a function near one of its roots is typically large and
this causes difficulties for the adaptive solver of   Section~\ref{section:solver}.

As with the phase function $\alpha_\nu$ and its derivative, the function
$\nu + \log(-Y_\nu(t)\sqrt{t})$ is represented via its values
at the $31$-point Chebyshev grid on a collection of subintervals
of (\ref{phase:50}).  It can be evaluated via barycentric Chebyshev
interpolation at any point on that interval, and the values
of $Y_\nu$ can obviously be obtained  from those of 
$\nu + \log(-Y_\nu(t)\sqrt{t})$.

\vskip 1em
{\it Stage three: computation of $-\nu + \log(J_\nu(t)\sqrt{t})$}

Assuming that $\nu > \frac{1}{2}$, we calculate the function
$-\nu + \log(J_\nu(t)\sqrt{t})$ on the interval
\begin{equation}
\left(0,\sqrt{\nu^2-\frac{1}{4}}\right].
\label{phase:200}
\end{equation}
This function is a solution
of the Riccati equation (\ref{phase:200}) 
with $q$  as in (\ref{phase:q}), and it is tempting to try
to calculate in the same way that  $\nu + \log(Y_\nu(t)\sqrt{t})$ 
is constructed in the preceding step.
That is, by evaluating $J_\nu$ and its derivative at a suitably
chosen point
\begin{equation}
t^{**} > \sqrt{\nu^2-\frac{1}{4}}
\end{equation}
and solving the corresponding terminal value problem for Riccati's
equation.  Such an approach is not numerically viable.  
The solution
\begin{equation}
-\nu + \log(J_\nu(t)\sqrt{t})
\label{phase:j}
\end{equation}
is recessive when solving the Riccati equation in the backward direction
while 
\begin{equation}
\nu + \log(-Y_\nu(t)\sqrt{t}) 
\label{phase:y}
\end{equation}
is dominant.  As a consequence, approximations of 
(\ref{phase:j}) obtained  by solving a terminal boundary value
problem for  (\ref{phase:riccati}) are highly inaccurate
while approximations of (\ref{phase:y}) obtained in such a fashion
are not.  
See, for instance, Chapter~I of \cite{GilSeguraTemme} for a discussion
of the recessive and dominant solutions of ordinary differential equations.

Rather than solving a terminal boundary value for (\ref{phase:riccati})
in order to calculate (\ref{phase:j}), we solve an initial value problem.
When $\nu \geq 10$, 
 we use the logarithm form (\ref{preliminaries:debye:expansion3}) of Debye's 
asymptotic of $J_\nu$ in order to evaluate 
(\ref{phase:j}) and its derivative at the left-hand endpoint
of (\ref{phase:200}).  When $\nu < 10$, Debye's expansion
is not necessarily sufficiently accurate and we use
the series expansion (\ref{preliminaries:series:4}) in order to evaluate
(\ref{phase:j}) and its derivative at the left-hand endpoint
of (\ref{phase:200}).
Again, our motivation for calculating
$-\nu + \log(J_\nu(t)\sqrt{t})$
in lieu of $\log(J_\nu(t)\sqrt{t})$ is that the former is bounded away from 
$0$ on the interval (\ref{phase:200}) while the latter is not.

The initial value problem is solved using the procedure of
Section~\ref{section:solver}, and,  as in the cases of 
$\alpha_\nu$ and  $\nu + \log(-Y_\nu(t)\sqrt{t})$,
we represent $-\nu + \log(J_\nu(t)\sqrt{t})$
via its value at the $31$-point Chebyshev grid on a collection of subintervals
of (\ref{phase:200}).
Using this data, the value of $J_\nu$ can be evaluated at any point in the interval $[a,b]$
via the obvious procedure.

\vskip 1em
\begin{remark}
Although the algorithm described in this section is highly specialized
to the case of Bessel's differential equation, it can, in fact, be modified
so as to apply to a large class of second order equations
of the form
\begin{equation}
y''(t) + q(t) y(t) = 0 \ \ \mbox{for all} \  \ a < t <b .
\label{phase:1000}
\end{equation}
Suppose, for instance, that  $q$ is smooth on $[a,b]$, has  a zero at $t_0 \in (a,b)$, 
is negative on  $(a,t_0)$  and is positive on  $(t_0,b)$. 
The procedure of the first stage for constructing
a nonoscillatory phase function on $(t_0,b)$ relies on an asymptotic expansion
which allows for the evaluation
of a nonoscillatory phase function at the point $b$.
In the absence of such an approximation, the algorithm of
\cite{BremerKummer} can be used instead.  That algorithm also
proceeds by solving Kummer's equation (\ref{preliminaries:riccati:kummer}),
but it incorporates a mechanism for numerically calculating
the appropriate initial values of a nonoscillatory phase function and its
derivatives.

The procedure of the second stage does not rely on any asymptotic or series
expansions of Bessel functions, only on the values of the phase function
computed in the first phase.  Consequently, it does not need to be
modified in order to obtain a solution of Riccati's equation
 which is increasing as $t \to 0^+$.

In the third stage, one of Debye's asymptotic expansions is
used to compute the values of the Bessel function $J_\nu$
and its derivative at a point near $0$.  In the event that such an approximation is 
not available, a solution of the Riccati equation which is increasing
as $t \to t_0$ from the left can be obtained by solving an initial
value problem with arbitrary initial conditions
and then scaling the result in order to make it consistent with 
the desired solution of (\ref{phase:1000}).
This procedure is analogous to that used in order to obtain
a recessive solution of a linear recurrence relation by running
the recurrence relation backwards (see, for instance, Section 3.6 of \cite{DLMF}).

Further generalization to the case in which $q$ has multiple zeros
on the interval $[a,b]$ is also possible, but beyond the scope
of this article.

\end{remark}

\label{section:phase}
\end{section}

\begin{section}{The numerical construction of the precomputed table}


In this section, we describe the  procedure used to construct
the table which allows for the numerical evaluation
of the Bessel functions $J_\nu$ and $Y_\nu$ for a large range
of parameters and arguments.
This table stores the coefficients in the piecewise
compressed bivariate Chebyshev expansions (as defined in Section~\ref{preliminaries:chebyshev4})
of several functions.
The functions $A_1$ and $C_1$  allow for
the evaluation of the nonoscillatory phase function $\alpha_\nu(t)$ 
 defined in Section~\ref{preliminaries:nonoscillatory}, as well as its
derivative $\alpha_\nu'(t)$, on the  subset
\begin{equation}
\mathcal{O}_1 = 
\left\{
(\nu,t) :
2 \leq \nu \leq 1\sep,000\sep,000\sep,000\ \ \mbox{and} \ \sqrt{\nu^2-\frac{1}{4}} \leq t \leq 1\sep,000 \nu
\right\}
\end{equation}
of the oscillatory region $\mathcal{O}$.
A second set of functions $A_2$ and $C_2$ allow for the evaluation
of $\alpha_\nu(t)$ and $\alpha_\nu'(t)$ on
\begin{equation}
\mathcal{O}_2 = 
\left\{
(\nu,t) :
0 \leq \nu \leq  2 \ \ \mbox{and}\ \ 2 \leq t \leq 1\sep,000
\right\}.
\end{equation}
A third set of functions  $B_1$ and $B_2$ allow for the evaluation  of $-\nu + \log(J_\nu(t)\sqrt{t})$
and $\nu + \log(-Y_\nu(t)\sqrt{t})$ on the subset 
\begin{equation}
\mathcal{N}_1 = 
\left\{
(\nu,t) :
\nu \geq 2   \ \ \mbox{and}\ \ 
\frac{\nu}{1\sep,000} < t < \sqrt{\nu^2-\frac{1}{4}}
\right\}
\label{phase:nonoscillatory}
\end{equation}
of the nonoscillatory region $\mathcal{N}$.
When $\nu$ is large, it is  numerically advantageous
to expand $\alpha_\nu$, $\alpha_\nu'$,
$-\nu + \log(J_\nu(t)\sqrt{t})$
and $\nu + \log(-Y_\nu(t)\sqrt{t})$ in powers of $\frac{1}{\nu}$ rather
than in powers of $\nu$.  Consequently, in this procedure the functions
$A_1$, $C_1$, $B_1$ and $B_2$ depend on  $x = \frac{1}{\nu}$.
Here, we only describe the construction of the functions $A_1$, $C_1$,
$B_1$ and $B_2$.  The procedure for the construction of $A_2$ and $C_2$
is quite similar, however.

There computations were conducted 
using IEEE extended precision
arithmetic in order to ensure high accuracy.    
The resulting table,
which consists of the coefficients in the expansions
of $A_1$, $C_1$, $A_2$, $C_2$, $B_1$ and $B_2$,
is approximately $1.3$ megabytes in size.
It allows for the evaluation of $\alpha_\nu$, $\alpha_\nu'$, $-\nu + \log(J_\nu(t)\sqrt{t})$
and $\nu + \log(-Y_\nu(t)\sqrt{t})$ with roughly double precision relative
accuracy (see the experiments of Section~\ref{section:experiments}).
The code was written in Fortran using OpenMP extensions and compiled with version 4.8.4 
of the GNU Fortran compiler.  It was executed on a computer
equipped with  $28$ Intel Xeon E5-$2697$ processor cores running at $2.6$ GHz.
The construction of this table took approximately $227$ seconds on this machine.

We conducted these calculations using extended precision arithmetic in order
to ensure that the resulting expansions obtained full double precision accuracy.
When these calculations are conducted in IEEE double precision arithmetic instead,
only a small amount of precision is lost.  We found that a table which
can evaluate $\alpha_\nu$, $\alpha_\nu'$, $-\nu + \log(J_\nu(t)\sqrt{t})$
and $\nu + \log(-Y_\nu(t)\sqrt{t})$
with roughly $12$ digits of relative accuracy could be constructed
using double precision arithmetic.  Less than $5$ seconds were required
to do so.


\vskip 1em
{\it Stage one: construction of the phase functions and logarithms}

\begin{table}
\def\arraystretch{1.3}
\center
\begin{tabular}{l@{\hspace{2em}}cc@{\hspace{3em}}l@{\hspace{2em}}cc}
\toprule
$j$ & $\xi_j$ & $\xi_{j+1}$ & $j$ & $\xi_j$ & $\xi_{j+1}$ \\
\midrule
\addlinespace[.25em]
0 &  $\frac{1}{1\sep,000\sep,000\sep,000}$ &  $ \frac{1}{100\sep,000\sep,000}$ 
&
5 &  $\frac{1}{10\sep,000}$ &  $ \frac{1}{1\sep,000}$ \\
1 &  $\frac{1}{100\sep,000\sep,000}$ &  $ \frac{1}{10\sep,000\sep,000}$  &
6 &  $\frac{1}{1\sep,000}$ &  $ \frac{1}{100}$ \\
2 &  $\frac{1}{10\sep,000\sep,000}$ &  $ \frac{1}{1\sep,000\sep,000}$ &
7 &  $\frac{1}{100}$ &  $ \frac{1}{50}$ \\
3 &  $\frac{1}{1\sep,000\sep,000}$ &  $ \frac{1}{100\sep,000}$ &
8 &  $\frac{1}{50}$ &  $ \frac{1}{10}$ \\
4 &  $\frac{1}{100\sep,000}$ &  $ \frac{1}{10\sep,000}$ &
9 &  $\frac{1}{10}$ &  $ \frac{1}{2}$ \\
\addlinespace[.25em]
\bottomrule
\end{tabular}
\caption{The endpoints of the intervals $\left[\xi_j,\xi_{j+1}\right]$
used in Stage one of the procedure of Section~\ref{section:expansions}.}
\label{expansion:table1}
\end{table}

We began this stage of the procedure by constructing a partition
\begin{equation}
\xi_0 < \xi_1 < \xi_2 < \ldots < \xi_{10} 
\label{expansion:100}
\end{equation}
of the interval
\begin{equation}
\left[ \frac{1}{1\sep,000\sep,000\sep,000}, \frac{1}{2}\right].
\label{expansion:101}
\end{equation}
This partition divides (\ref{expansion:101}) into
ten subintervals, the endpoints of which are given in Table~\ref{expansion:table1}.
For each such interval $\left[\xi_j,\xi_{j+1}\right]$, we formed
the nodes
\begin{equation}
x_1^{(j)}, \ldots, x_{50}^{(j)}
\label{expansion:102}
\end{equation}
of the $50$-point Chebyshev  grid on $\left[\xi_j,\xi_{j+1}\right]$.
Next, for each $x$ in the collection
\begin{equation}
x_1^{(0)}, \ldots, x_{50}^{(0)},x_1^{(1)}, \ldots, x_{50}^{(1)}, \ldots,
x_1^{(9)}, \ldots, x_{50}^{(9)}
\label{expansion:103}
\end{equation}
we executed
the algorithm of Section~\ref{section:phase} with $\nu$ take to be $\frac{1}{x}$.
The requested precision for the solver of Section~\ref{section:solver}
used by the algorithm of Section~\ref{section:phase}
was set to $\epsilon = 10^{-20}$ and we set the parameter $n$ to be $50$
so that the functions produced by the algorithm of Section \ref{section:solver}
are represented via their values on the $50$-point Chebyshev grids on 
a collection of subintervals.    Were it not for the fact that 
the solver of (\ref{section:solver})
runs in time independent of $\nu$, these calculations would be prohibitely
expensive to carry out, even on a massively parallel computer.

For each value of $\nu$ corresponding to one of the points
(\ref{expansion:103}), this results in the values 
of $\alpha_\nu$ and  $\alpha_\nu'$ at the nodes
of the $50$-point Chebyshev grids on a collection of subintervals of 
\begin{equation}
\left[\sqrt{\nu^2-\frac{1}{4}},1000 \nu\right)
\label{expansion:104}
\end{equation}
and piecewise Chebyshev expansions for $-\nu + \log(J_\nu(t)\sqrt{t})$
and $\nu + \log(-Y_\nu(t)\sqrt{t})$ on a collection of subintervals
of 
\begin{equation}
\left[\frac{\nu}{1000},
\sqrt{\nu^2-\frac{1}{4}}\right].
\label{expansion:105}
\end{equation}
Using this data and the techniques
discussed in Section~\ref{preliminaries:chebyshev3}, 
we can evaluate $\alpha_\nu$ and its derivative at any point
in (\ref{expansion:104}) and we can evaluate
 $-\nu + \log(J_\nu(t)\sqrt{t})$
and $\nu + \log(-Y_\nu(t)\sqrt{t})$ at any point in (\ref{expansion:105}).


\vskip 1em
{\it Stage two: formation of unified discretizations}

Now, for each $x$ in the set (\ref{expansion:103}) we adaptively discretize
the function $f_x:[0,1] \to \mathbb{R}$ defined via
\begin{equation}
f_x(y) = \alpha_\nu \left(\sqrt{\nu^2-\frac{1}{4}} + \left(1\sep,000\nu - \sqrt{\nu^2-\frac{1}{4}}\right)y \right)
\ \ \mbox{with}\ \ \nu = \frac{1}{x}
\label{expansion:106}
\end{equation}
using the procedure of  Section~\ref{preliminaries:adaptive}.  We ask for
$\epsilon = 10^{-17}$ precision and take the parameter $n$ to be $49$.  Each discretization
consists of a collection of subintervals of $[0,1]$ on which $f_x$ can be represented
to high accuracy using a $49$-term Chebyshev expansion.  We form a unified
discretization
\begin{equation}
\left[a_0,a_1\right], \left[a_1,a_2\right], \ldots, \left[a_{24},a_{25}\right]
\label{expansion:107}
\end{equation}
of $[0,1]$ by merging these discretizations; that is, by ensuring that 
the sets (\ref{expansion:107}) have the property that
each subset appearing in the discretization of one of the functions $f_x$
is the union of some collection of the subintervals (\ref{expansion:107}).

For each $x$, we also adaptively discretize each of the functions
$g_x:[0,1] \to \mathbb{R}$ and $h_x:[0,1]\to\mathbb{R}$ 
defined via the formulas
\begin{equation}
g_x(y) = 
-\nu + \log\left(J_\nu\left(t\right)\sqrt{t}\right)
\ \ \mbox{with}\ \ \nu = \frac{1}{x},
 \ t = \frac{\nu}{1\sep,000} + 
\left(\sqrt{\nu^2-\frac{1}{4}} - \frac{\nu}{1\sep,000}\right) y
\label{expansion:108}
\end{equation}
and
\begin{equation}
h_x(y) = 
\nu + \log\left(-Y_\nu\left(t\right)\sqrt{t}\right)
\ \ \mbox{with}\ \ \nu = \frac{1}{x},
 \ t = \frac{\nu}{1\sep,000} + 
\left(\sqrt{\nu^2-\frac{1}{4}} - \frac{\nu}{1\sep,000}\right) y.
\label{expansion:109}
\end{equation}
Again, we used the procedure of Section~\ref{preliminaries:adaptive}
with  
$\epsilon = 10^{-17}$ and  $n = 49$. 
We then formed a unified discretization
\begin{equation}
\left[b_0,b_1\right], \left[b_1,b_2\right], \ldots, \left[b_{22},b_{23}\right]
\label{expansion:110}
\end{equation}
of $[0,1]$ in the same fashion in which we formed (\ref{expansion:107}).

\vskip 1em
{\it  Stage three: construction of the functions $A_1$ and $C_1$}

We define $A_1$ via the formula
\begin{equation}
A_1(x,y) = \frac{1}{\nu} 
\alpha_\nu \left(\sqrt{\nu^2-\frac{1}{4}} + \left(1\sep,000\nu - \sqrt{\nu^2-\frac{1}{4}}\right)y \right)
\ \ \mbox{with}\ \ \nu = \frac{1}{x}
\end{equation}
and $C_1$ via the formula
\begin{equation}
C_1(x,y) = \frac{1}{\nu} \alpha_\nu' 
\left(\sqrt{\nu^2-\frac{1}{4}} + \left(1\sep,000\nu - \sqrt{\nu^2-\frac{1}{4}}\right)t \right)
\ \ \mbox{with}\ \ \nu = \frac{1}{x}.
\end{equation}
Obviously, $A_1$ and $C_1$ are defined on the compact rectangle
\begin{equation}
\left[ \frac{1}{1\sep,000\sep,000\sep,000}, \frac{1}{2}\right]
\times [0,1].
\end{equation}
For each $i=0,\ldots,9$ and each $j=0,\ldots,24$, we form the 
$49^{th}$ order  compressed bivariate Chebyshev expansions of $A_1$  and $C_1$ on the
compact rectangle
\begin{equation}
\left[\xi_i,\xi_{i+1}\right] \times \left[a_j,a_{j+1} \right].
\end{equation}
There are $250$ such rectangles and the uncompressed bivariate Chebyshev expansions
of order $49$ on each rectangle would involve $2\sep,500$ coefficients.   A total
of  $250 \times 2\sep,500 = 625\sep,000$ coefficients
would be required to store the uncompressed bivariate expansions
for $A_1$, and another  $625\sep,000$ would be required for $C_1$.
The compressed bivariate expansions are much smaller.  A mere
 $31\sep,884$ values (this includes both the coefficients and 
the indices appearing in the sums (\ref{preliminaries:chebyshev4:expansion}), which must also be stored)
were  required to represent $A_1$.  Only  $51\sep,076$ values were 
needed to represent $C_1$.

\vskip 1em
{\it  Stage four: construction of the functions $B_1$ and $B_2$}

We define $B_1$ via the formula
\begin{equation}
B_1(x,y) = 
-1 + \frac{1}{\nu}\log\left(J_\nu\left(t\right)\sqrt{t}\right)
\ \ \mbox{with}\ \ \nu = \frac{1}{x},
 \ t = 
\frac{\nu}{1\sep,000} + \left(\sqrt{\nu^2-\frac{1}{4}}-\frac{\nu}{1\sep,000}\right)y 
\end{equation}
and $B_2$ via the formula
\begin{equation}
B_1(x,y) = 
1 + \frac{1}{\nu}\log\left(-Y_\nu\left(t\right)\sqrt{t}\right)
\ \ \mbox{with}\ \ \nu = \frac{1}{x},
 \ t = 
\frac{\nu}{1\sep,000} + \left(\sqrt{\nu^2-\frac{1}{4}}-\frac{\nu}{1\sep,000}\right)y .
\end{equation}
Obviously, $B_1$ and $B_2$ are defined on the compact rectangle
\begin{equation}
\left[ \frac{1}{1\sep,000\sep,000\sep,000}, \frac{1}{2}\right]
\times [0,1].
\end{equation}
For each $i=0,\ldots,9$ and each $j=0,\ldots,23$, we form the 
$49^{th}$ order  compressed bivariate Chebyshev expansions of $B_1$  and $B_2$ on the
compact rectangle
\begin{equation}
\left[\xi_i,\xi_{i+1}\right] \times \left[b_j,b_{j+1} \right].
\end{equation}
There are $230$ such rectangles and the uncompressed bivariate Chebyshev expansions
of $B_1$ and $B_2$ would invole $2 \times 230 \times 2\sep,500 = 1\sep,150\sep,000$ coefficients.
Using the compressed expansions, we are able to store $B_1$ using
$32\sep,910$ values and $B_2$ with $46\sep,950$ values.

\vskip 1em
\begin{remark}
It is possible to compute the values of both $\alpha_\nu$ and $\alpha_\nu'$
using the function $A_1$ via spectral differentiation   (as discussed, 
for instance, \cite{Mason-Hanscomb}). Such an approach would, however,
lead a level of  loss of precision in the obtained values of $\alpha_\nu'$
which we find unacceptable.  

In a similar vein, spectral integration
 could be used to evaluate
$\alpha_\nu$ given the values of, or an expansion for, $\alpha_\nu'$.  
Spectral integration does not suffer from the same defect as 
spectral differentiation and such a calculation could be carried out
with little loss of precision;
however, integration of $\alpha_\nu'$ 
can only determine  $\alpha_\nu$   up to a constant.  
The appropriate constant would have to be calculated or stored
in some fashion.  We chose the simpler, but possibly more expensive,
procedure described in this paper over such approach.
\end{remark}

In order to evaluate $\alpha_\nu(t)$ and $\alpha_\nu'(t)$ given the expansions
of $A_1$ and $C_1$ constructed using the procedure describe above,
we execute the following sequence of steps:

\begin{enumerate}

\item
First, we  let 
\begin{equation}
x = \frac{1}{\nu}
\end{equation}
and
\begin{equation}
y = \frac{t - \sqrt{\nu^2-\frac{1}{4}}}{  \left(1\sep,000\nu - \sqrt{\nu^2-\frac{1}{4}}\right)}.
\end{equation}

\vskip 1em
\item
Next, we next search through the intervals (\ref{expansion:103}) 
in order to find the index $i$ of the one containing $x$ and through the intervals
(\ref{expansion:107}) for index $j$ of the interval containing $y$. 

\vskip 1em
\item
Having discovered that $(x,y) \in  \left[\xi_i,\xi_{i+1}\right] \times \left[a_j,a_{j+1}\right]$,
we evaluate the compressed bivariate Chebyshev series expansion representing
$A_1$ on this rectangle.  We scale the result by $\nu$ in order to
obtain the value of $\alpha_\nu(t)$.  We then evaluate the compressed
bivariate Chebyshev expansion representing $C_1$ on this rectangle.  We scale
the result by $\nu$ in order to obtain the value of $\alpha_\nu'(t)$.

\end{enumerate}
A virtually identical procedure is used to evaluate $\log(J_\nu(t))$ and
$\log(-Y_\nu(t))$ using the expansions of $B_1$ and $B_2$ stored in the 
table.

\label{section:expansions}
\end{section}

\begin{section}{An  algorithm for the rapid numerical evaluation of Bessel functions}

In this section, we describe the operation of our code for evaluating
the Bessel functions $J_\nu$ and $Y_\nu$ of nonnegative orders and positive
arguments.  It was written in Fortran and its interface to the user consists of two
subroutines, one called {\tt bessel\_eval_init} and the other {\tt bessel\_eval}.
The {\tt bessel\_eval\_init} routine reads the precomputed table constructed
via the procedure of Section~\ref{section:expansions}  from the disk into memory.
Once this has been accomplished, the {\tt bessel\_eval} can be called.  
It takes as input an order $\nu \geq 0$ and  an argument $t >0$.
When $(\nu,t)$ is in the oscillatory region $\mathcal{O}$, 
it returns the values of $\alpha_\nu(t)$ and $\alpha'(t)$ as well as
those of $J_\nu(t)$ and $Y_\nu(t)$.  When $(\nu,t)$ is in the 
nonoscillatory region $\mathcal{N}$, it returns the values
of $\log(J_\nu(t))$ and $\log(-Y_\nu(t))$ as well as those
of $J_\nu(t)$ and $Y_\nu(t)$.  Of course, when $t \ll \nu$,
these latter values might not be representable via the IEEE double format
arithmetic.  In this event, $0$ is returned for $J_\nu(t)$ and
$-\infty$ for $Y_\nu(t)$.

The {\tt bessel\_eval} code as well as the  code for all of the experiments described
in the following section are available from the GitHub repository at address

\hskip 2em {\tt http://github.com/JamesCBremerJr/BesselEval}

and from the author's website at the address

\hskip 2em {\tt http://www.math.ucdavis.edu/\textasciitilde bremer/code.html}.

The code for the experiments described in  Sections~\ref{experiments:4}, \ref{experiments:5}
and \ref{experiments:6} depends on the package described in \cite{Amos}.  
Fortran source code for \cite{Amos} can be obtained from Netlib's  TOMS repository:

\hskip 2em {\tt http://netlib.org/toms/}.

\vskip 1em

The {\tt bessel\_eval} code operates as follows:
\begin{enumerate}

\item
When $\nu \geq 2$ and $\sqrt{\nu^2-\frac{1}{4}} \leq t \leq 1\sep,000 \nu$, the precomputed
expansions of $A_1$ and its derivative are used to evaluate the nonoscillatory
phase function $\alpha_\nu$ and its derivative $\alpha_\nu'$ at the point $t$.
Then, formulas (\ref{preliminaries:nonoscillatory:1}) and 
(\ref{preliminaries:nonoscillatory:2})  are used to produce the values
of $J_\nu(t)$ and $Y_\nu(t)$.

\vskip 1em
\item
When $\nu \geq 2$ and $\frac{\nu}{1\sep,000} \leq  t < \sqrt{\nu^2-\frac{1}{4}}$,  precomputed
expansions of $B_1$ and $B_2$ are used to evaluate 
$-\nu + \log(J_\nu(t)\sqrt{t})$  and $\nu + \log(Y_\nu(t)\sqrt{t})$.  The values
of $J_\nu(t)$ and $Y_\nu(t)$ are calculated in the obvious fashion.
Note that it is the the values of  $\log(J_\nu(t))$ and $\log(-Y_\nu(t))$
and  not those of $-\nu + \log(J_\nu(t)\sqrt{t})$  and $\nu + \log(Y_\nu(t)\sqrt{t})$
that are returned by the  {\tt bessel\_eval} routine.  
\vskip 1em
\item
When $\nu \geq 2$ and $t \leq \frac{\nu}{1\sep,000}$, 
Debye's expansions (\ref{preliminaries:series:4})
and (\ref{preliminaries:series:5})
are used to evaluate $\log(J_\nu(t))$ and $\log(-Y_\nu(t))$.
The values of  $J_\nu(t)$ and $-Y_\nu(t)$ are computed as one would expect.


\vskip 1em
\item
When $\nu < 2$ and $2 \leq t \leq 1\sep,000$, the precomputed
expansions of $A_2$ and $C_2$ are used to evaluate the nonoscillatory
phase function $\alpha_\nu$ and its derivative $\alpha_\nu'$ at the point $t$.
Then, formulas (\ref{preliminaries:nonoscillatory:1}) and 
(\ref{preliminaries:nonoscillatory:2})  are used to produce the values
of $J_\nu(t)$ and $Y_\nu(t)$.

\vskip 1em
\item
When $\nu < 2$, $t < 2$ and  $(\nu,t)$ is in 
the oscillatory region,
 we use the series expansions (\ref{preliminaries:series:1})
and (\ref{preliminaries:series:2}) in 
 in order to evaluate $J_\nu(t)$ and $Y_\nu(t)$.
As discussed in Section~\ref{preliminaries:series},
Chebyshev interpolation is used in the computation
of $Y_\nu$ when $\nu$ is either an integer or
close to one.  The value of $\alpha_\nu'$
is evaluated via the formula (\ref{preliminaries:nonoscillatory:alphap}),
and that of $\alpha_\nu$ is calculated via
\begin{equation}
\alpha_\nu(t) = \arctan\left(\frac{Y_\nu(t)}{J_\nu(t)}\right).
\end{equation}

\vskip 1em
\item
When $\nu < 2$, $t < 2$ and  $(\nu,t)$ is in 
the nonoscillatory regime,
 the series expansions (\ref{preliminaries:series:3})
and (\ref{preliminaries:series:4}) are used to 
evaluate $\log(J_\nu(t))$ and $\log(-Y_\nu(t))$.
As discussed in Section~\ref{preliminaries:series},
Chebyshev interpolation is used in the computation
of $Y_\nu$ when $\nu$ is either an integer or
close to one.  The values of $J_\nu(t)$ and $Y_\nu(t)$
are calculated in the obvious fashion.
We use series expansions rather than Debye's expansion to evaluate
$J_\nu(t)$ and $Y_\nu(t)$ in this case because Debye's expansions
lose accuracy when $\nu$ is small.

\end{enumerate}

\label{section:numerics}
\end{section}


\begin{section}{Numerical Experiments}

In this section, we describe the results of numerical experiments conducted to
illustrate the performance of the {\tt bessel\_eval} subroutine.
These experiments were carried out on a laptop computer equipped
with an Intel Core i7-5600U processor running at 2.6 GHz and 16GB of memory.
Our code was compiled with the GNU Fortran compiler version 5.2.1 using
the ``-O3'' compiler optimization flag.


\begin{subsection}{The accuracy with which  $\alpha_\nu'$ is evaluted in the oscillatory region}

In these experiments, we measured the accuracy with which {\tt bessel\_eval}
calculates values of  $\alpha_\nu'$  in the oscillatory region.
We did so by comparison with highly accurate reference values computed using version 11
of Wolfram's Mathematica package.  

In each experiment, we first constructed $1\sep,000$ pairs $(\nu,t)$ by first
choosing $\nu$ in a given range and then randomly selecting a value of $t$ in 
the interval  
$$\left( \sqrt{\nu^2-1/4}, 1000\ \nu \right).$$
Unless, that is,
 $\nu < 1/2$, in which case we selected a  random vaue of $t$ in the interval $(0,1000)$ instead.  
For each pair $(\nu,t)$ obtained in this fashion, we calculated the relative
error in the value of  $\alpha_\nu'(t)$ returned by  bessel\_eval.
Table~\ref{table1} displays the results.  There, each row corresponds
 to one experiment and reports the maximum observed relative error
as well as the average running time of {\tt bessel\_eval}.

\begin{table}[h!!!]
\center
\small
 \begin{tabular}{lcc}
 \toprule
 Range of $\nu$ & Maximum relative         & Average evaluation\\
                & error in $\alpha_\nu'(t)$   & time (in seconds) \\
 \midrule
0 - 1 & 1.99\e{-15} & 2.06\e{-07}  \\
1 - 10 & 4.44\e{-16} & 9.99\e{-08}  \\
10 - 100 & 1.11\e{-16} & 1.03\e{-07}  \\
100 - 1\sep,000 & 1.11\e{-16} & 1.15\e{-07}  \\
1\sep,000 - 10\sep,000 & 1.11\e{-16} & 1.13\e{-07}  \\
10\sep,000 - 100\sep,000 & 1.11\e{-16} & 1.13\e{-07}  \\
100\sep,000 - 1\sep,000\sep,000 & 1.11\e{-16} & 1.17\e{-07}  \\
1\sep,000\sep,000 - 10\sep,000\sep,000 & 1.11\e{-16} & 1.18\e{-07}  \\
10\sep,000\sep,000 - 100\sep,000\sep,000 & 1.11\e{-16} & 1.17\e{-07}  \\
100\sep,000\sep,000 - 1\sep,000\sep,000\sep,000 & 2.22\e{-16} & 1.62\e{-07}  \\
 \bottomrule
 \end{tabular}

 \caption{The results of the experiments of Section~\ref{experiments:1} in which
 the accuracy of the evaluation of $\alpha_\nu'$  in the oscillatory region
 is tested through comparison with highly accurate reference values.}
 \label{table1}
\end{table}

\label{experiments:1}
\end{subsection}


\begin{subsection}{The accuracy with which  $-\nu+\log(J_\nu(t))$ and 
$-\nu + \log(-Y_\nu(t))$ are evaluated for small to moderate values of $\nu$}

In these experiments, we measured the accuracy with which {\tt bessel\_eval}
calculates  $-\nu+\log(J_\nu(t))$ and  $-\nu + \log(-Y_\nu(t))$
  in the nonoscillatory region.  Reference values for these experiments were generated 
using version 11.0 of  Wolfram's Mathematica package.  A considerable
amount of  time is required
for Mathematica to evaluate the Bessel functions $J_\nu(t)$ and $Y_\nu(t)$
when the magnitude of $\nu$ is large and 
$t$ is small relative to $\nu$.
Consequently, in these experiments we only considered values of  
$\nu$ between $\frac{1}{2}$ and $10\sep,000$.      Larger values of $\nu$ were treated
in the experiments described in the following section.

\begin{table}[h!!!]
\small\center
 \begin{tabular}{lccc}
 \toprule
 Range of $\nu$ & Maximum relative    & Maximum relative & Average evaluation\\
                & error in            & error in & time (in seconds) \\
                &$-\nu + \log(J_\nu(t))$ & $\nu + \log(-Y_\nu(t))$   & \\
 \midrule
0.5
 - 1 & 4.11\e{-16} & 7.01\e{-15} & 9.40\e{-07}  \\
1 - 10 & 2.44\e{-15} & 8.51\e{-15} & 1.06\e{-06}  \\
10 - 100 & 2.01\e{-15} & 3.16\e{-15} & 6.99\e{-07}  \\
100 - 1\sep,000 & 3.59\e{-15} & 3.74\e{-15} & 6.14\e{-07}  \\
1\sep,000 - 10\sep,000 & 1.70\e{-15} & 2.64\e{-15} & 3.87\e{-07}  \\
 \bottomrule
 \end{tabular}

 \caption{The results of the experiments of Section~\ref{experiments:2} in which
 the accuracy of {\tt bessel\_eval} in the nonoscillatory region
 is tested via comparison with highly accurate reference values
generated using Wolfram's Mathematica package.}
\label{table2}
\end{table}

In each experiment, we constructed $1\sep,000$ pairs $(\nu,t)$
by first chosing a value of $\nu$ in a given range and then
selecting a random point $t$ in the interval
$$\left(0,\sqrt{\nu^2-1/4}\right).$$
%
For each pair $(\nu,t)$ obtained in this fashion, we calculated
the relative accuracy of the quantities 
$-\nu + \log(J_\nu(t))$ and  $\nu +  \log(-Y_\nu(t))$.
Table~\ref{table2} displays the results of these experiments.  There each row corresponds
to one experiment and reports the largest relative errors which were
observed as well as the average time taken by the {\tt bessel\_eval} routine.

\label{experiments:2}
\end{subsection}

\begin{subsection}{The accuracy of the evaluation of $-\nu + \log(J_\nu(t))$ and $\nu+\log(-Y_\nu(t))$ 
deep in the nonoscillatory region}

The {\tt bessel\_eval} subroutine makes use of the
 asymptotic expansions  (\ref{preliminaries:debye:expansion1})
and (\ref{preliminaries:debye:expansion2}) when $0 < t < \nu / 10\sep,000 $.
For large $\nu$, Debye expansion's are accurate in a much larger interval.
In these experiments, we exploit this fact in  order to measure the 
accuracy with which {\tt bessel\_eval} calculates $-\nu + \log(J_\nu(t))$ and
$\nu + \log(-Y_\nu(t))$ deep in the nonoscillatory region.

In each experiment, we constructed $1\sep,000$ pairs by first selecting
a value of $\nu$ in a given range and thne 
picking a random value of $t$ in the interval
$\left(\nu/1000,\nu/10\right)$.  For each pair $(\nu,t)$ obtained in this
fashion, we computed the values of both $-\nu+\log(J_\nu(t))$ and $\nu + \log(-Y_\nu(t))$
using {\tt bessel\_eval}
and   compared them to reference values obtained using Debye's expansion.
The reference calculations were performed using IEEE quadruple precision arithmetic
in order to ensure high accuracy.
The results are shown in Table~\ref{table3}.
Each row there corresponds to one experiment and reports the range of $\nu$,
the maximum relative error which was observed,
 and the average time taken by {\tt bessel\_eval}.

 \begin{table}[h!!]
\center
 \begin{tabular}{lccc}
 \toprule
 Range of $\nu$ & Maximum relative    & Maximum relative & Average evaluation\\
                & error in            & error in & time (in seconds) \\
                &$-\nu + \log(J_\nu(t))$ & $\nu + \log(-Y_\nu(t))$   & \\
 \midrule
100 - 1\sep,000 & 1.53\e{-15} & 1.44\e{-15} & 3.32\e{-07}  \\
1\sep,000 - 10\sep,000 & 1.21\e{-15} & 1.79\e{-15} & 2.91\e{-07}  \\
10\sep,000 - 100\sep,000 & 1.26\e{-15} & 1.23\e{-15} & 2.59\e{-07}  \\
100\sep,000 - 1\sep,000\sep,000 & 1.02\e{-15} & 1.00\e{-15} & 2.49\e{-07}  \\
1\sep,000\sep,000 - 10\sep,000\sep,000 & 7.38\e{-15} & 7.46\e{-15} & 7.11\e{-07}  \\
10\sep,000\sep,000 - 100\sep,000\sep,000 & 1.02\e{-15} & 1.20\e{-15} & 4.58\e{-07}  \\
100\sep,000\sep,000 - 1\sep,000\sep,000\sep,000 & 1.25\e{-15} & 1.01\e{-15} & 3.24\e{-07}  \\
 \bottomrule
 \end{tabular}

 \caption{The results of the experiments of Section~\ref{experiments:3} in which
 the accuracy with which $-\nu + \log(J_\nu(t))$ and  $\nu + \log(-Y_\nu(t))$  
is tested for values of $100 \leq \nu  \leq 1\sep,000\sep,000\sep,000$
and $t \ll \nu$ through comparison with Debye's expansions.}
\label{table3}
\end{table}

\label{experiments:3}
\end{subsection}

\begin{subsection}{The accuracy of the evaluation of $J_\nu(t)$ and $Y_\nu(t)$ as a function
of $t$}

In these experiments, 
we measured the relative accuracy with which
{\tt bessel\_eval} calculates the Hankel function of the first
kind $H_\nu(t) = J_\nu(t) + i Y_\nu(t)$ as
a function of $t$.  We considered the Hankel function
instead of treating $J_\nu$ and $Y_\nu$ separately because $H_\nu(t)$
 does not vanish in the interval $(0,\infty)$
and its absolute value is nonoscillatory there, properties
not shared by the Bessel functions $J_\nu$ and $Y_\nu$.

In each experiment, we chose a value of $\nu$ and 
measured the relative accuracy with which {\tt bessel\_eval} calculates
 $J_\nu(t) + iY_\nu(t)$ at each of $1\sep,000$ equispaced points in the interval 
$\left[\nu,100000 \nu \right]$.  
Highly accurate reference values for these experiments were computed using 
Mathematica.  We chose the following values of 
$\nu$: $\sqrt{2}$, $10\sqrt{2}$, $100\sqrt{2}$ and $1\sep,000\sqrt{2}$ .
We also repeated these experiments using Amos' well-known code \cite{Amos}.

Figure~\ref{figure1} displays the results.    Each graph there
plots the base-$10$ logarithms of the relative errors in the calculated
values of $J_\nu(t) + i Y_\nu(t)$ as dots.
The graph of the function  $\kappa(t) \epsilon_0$, where $\kappa(t)$ is the 
condition number of the evaluation of $H_\nu$
at the point $t$ and $\epsilon_0 = 2^{-52} \approx 2.22044604925031e{-16}$
 is machine epsilon, is also plotted as a solid curve.
The results for {\tt bessel\_eval} are shown on the left while those for
Amos' code appear on the right on the right.

\begin{figure}[h!!!]
\begin{center}
\includegraphics[width=.40\textwidth]{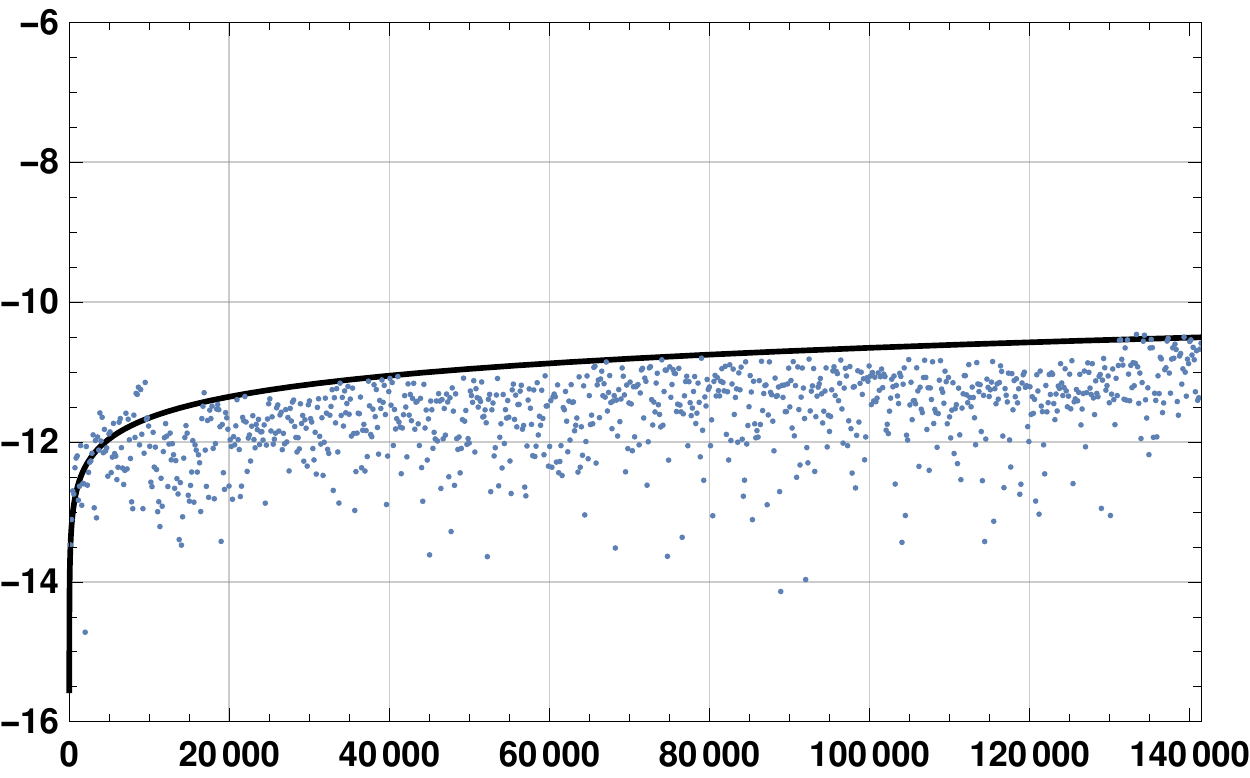}
\hfil
\includegraphics[width=.40\textwidth]{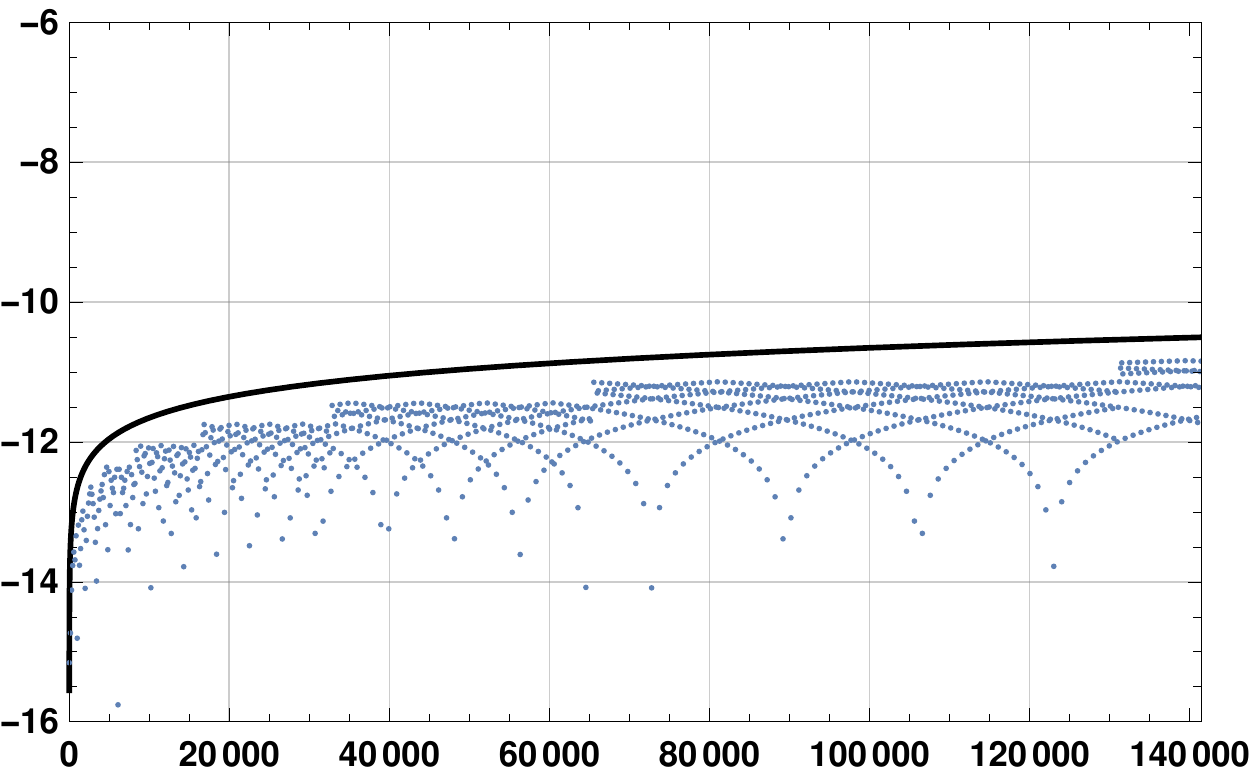}
{\center\small $\nu = \sqrt{2}$}
\vskip .7em

\includegraphics[width=.40\textwidth]{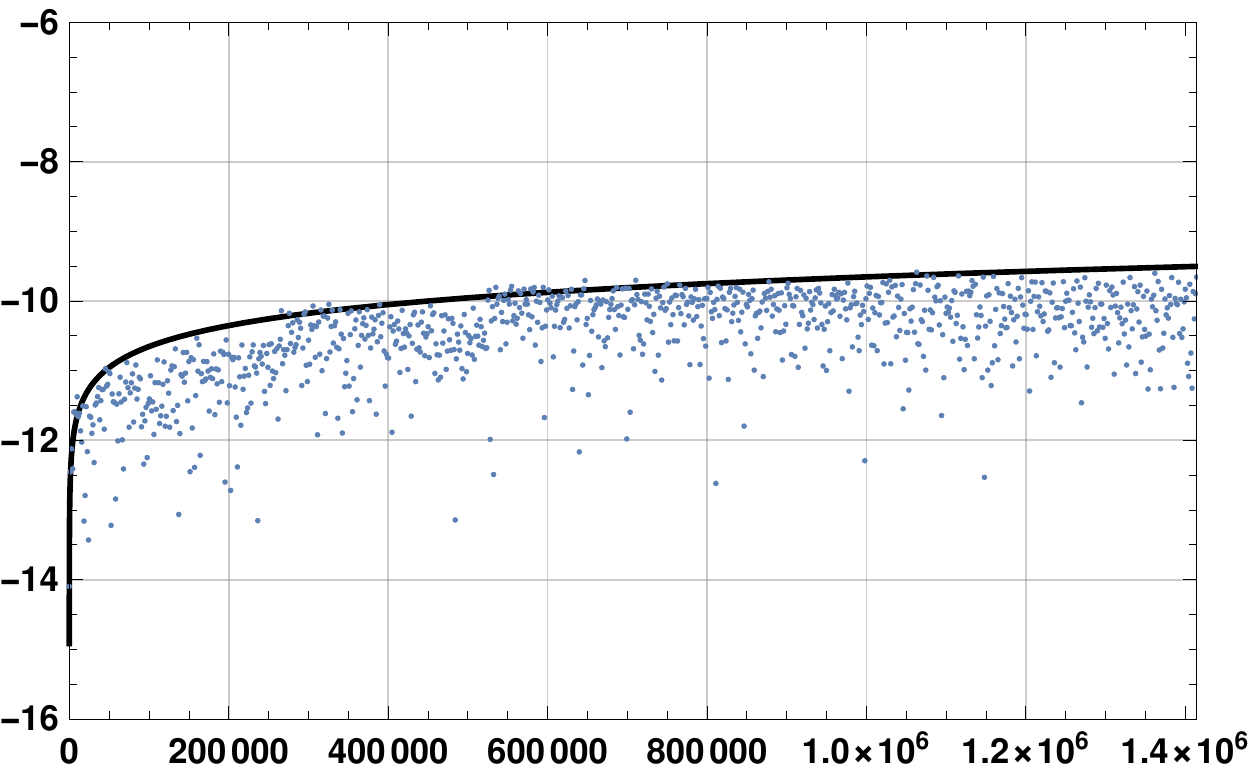}
\hfil
\includegraphics[width=.40\textwidth]{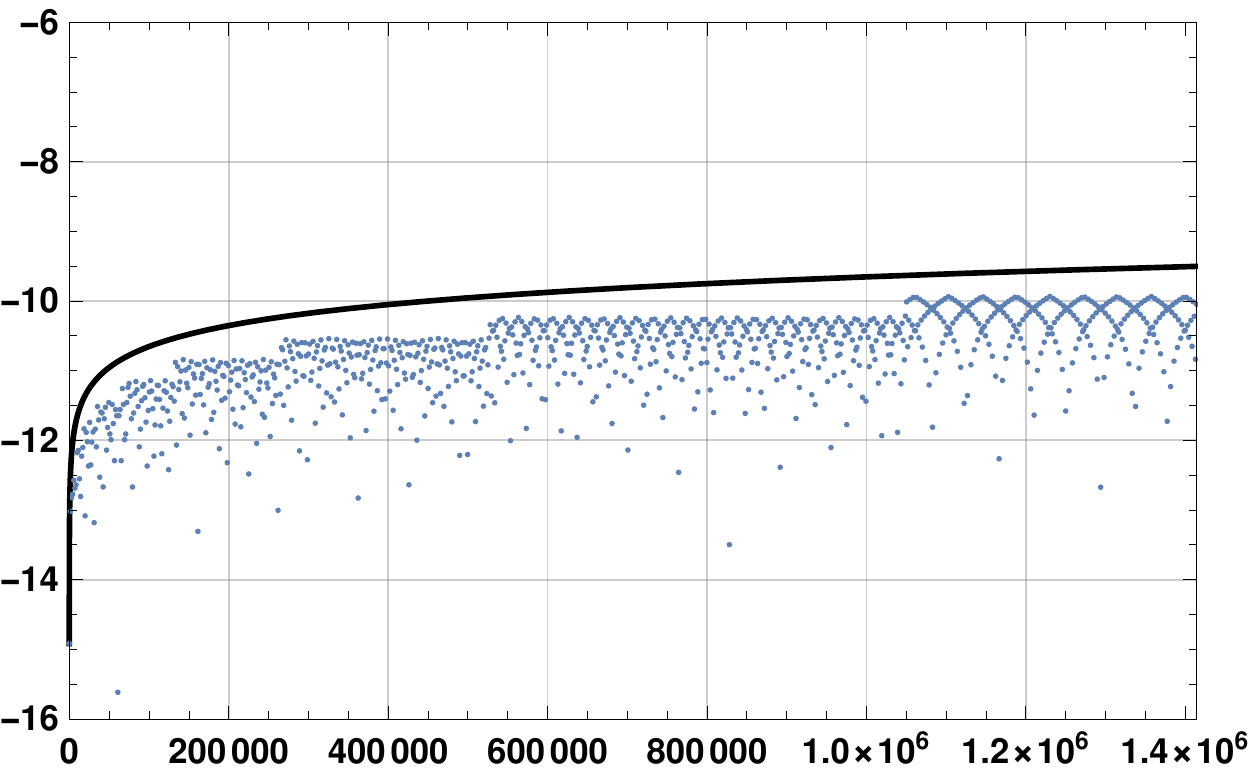}
{\center\small $\nu = 10\sqrt{2}$}
\vskip .7em

\includegraphics[width=.40\textwidth]{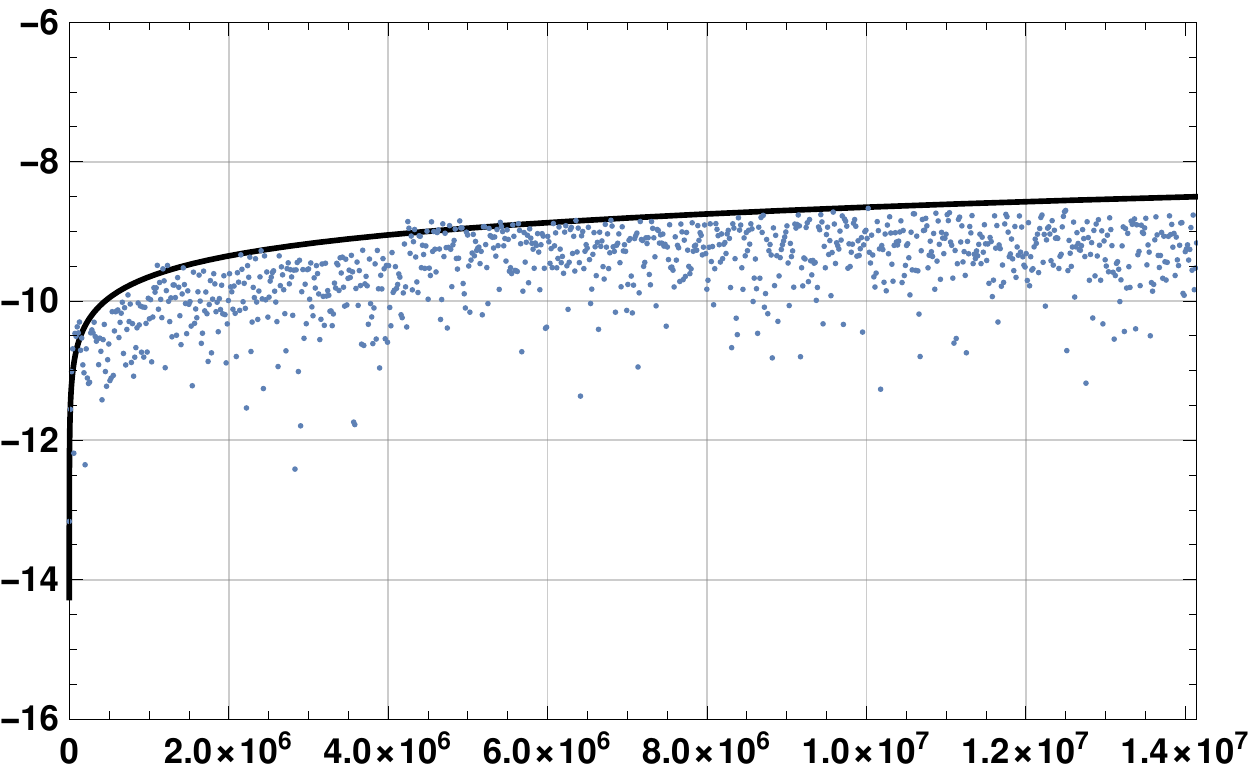}
\hfil
\includegraphics[width=.40\textwidth]{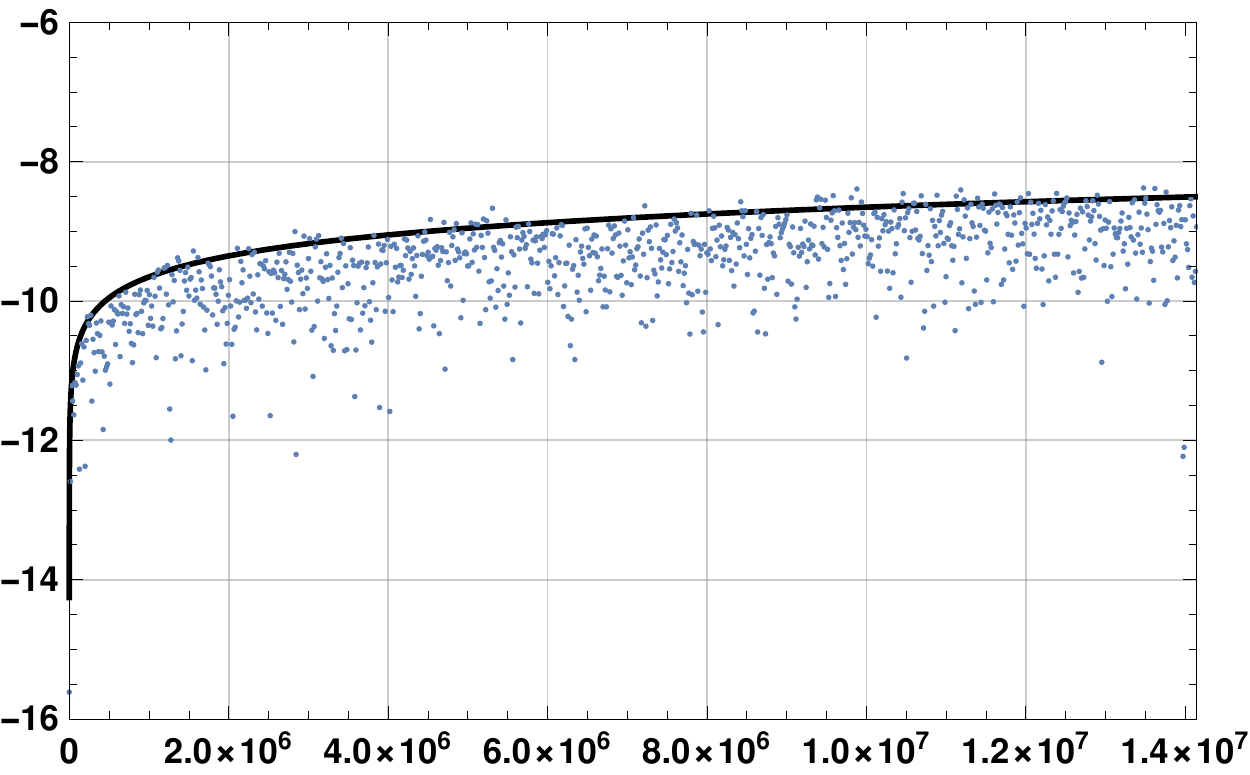}
{\center\small$\nu = 100\sqrt{2}$}
\vskip .7em

\includegraphics[width=.40\textwidth]{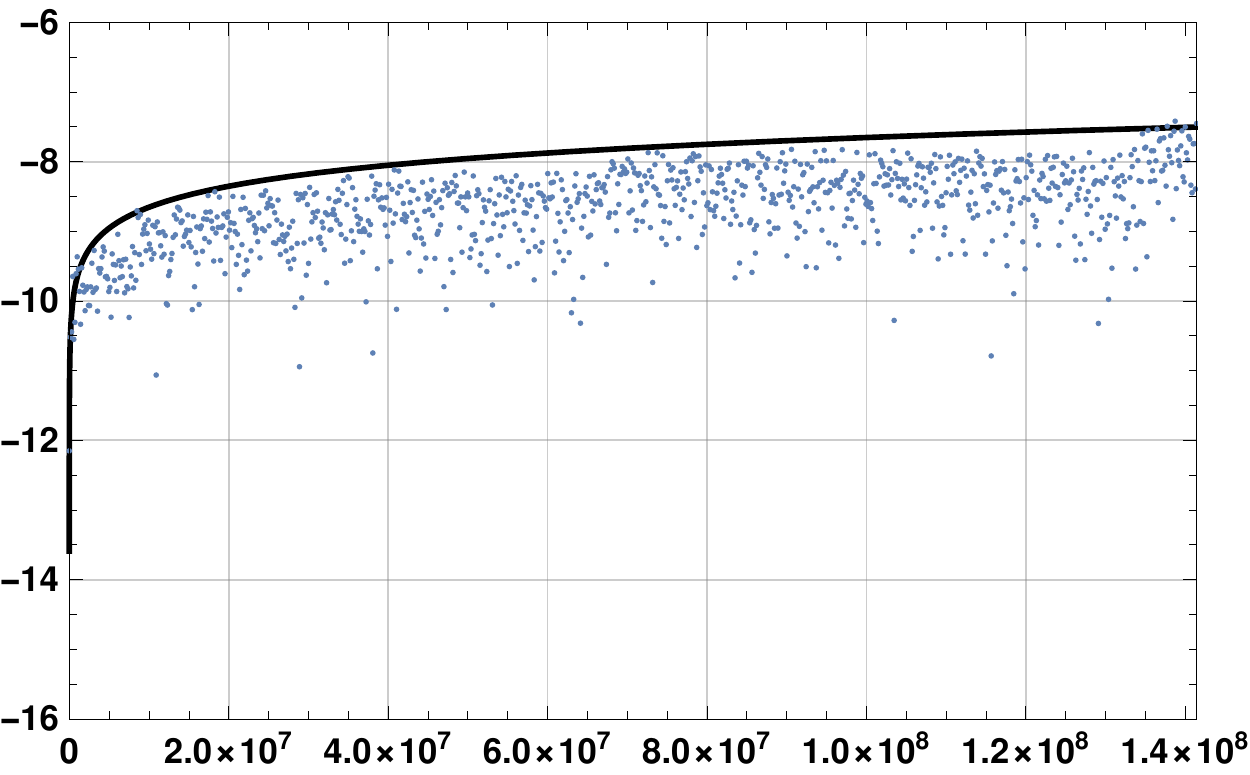}
\hfil
\includegraphics[width=.40\textwidth]{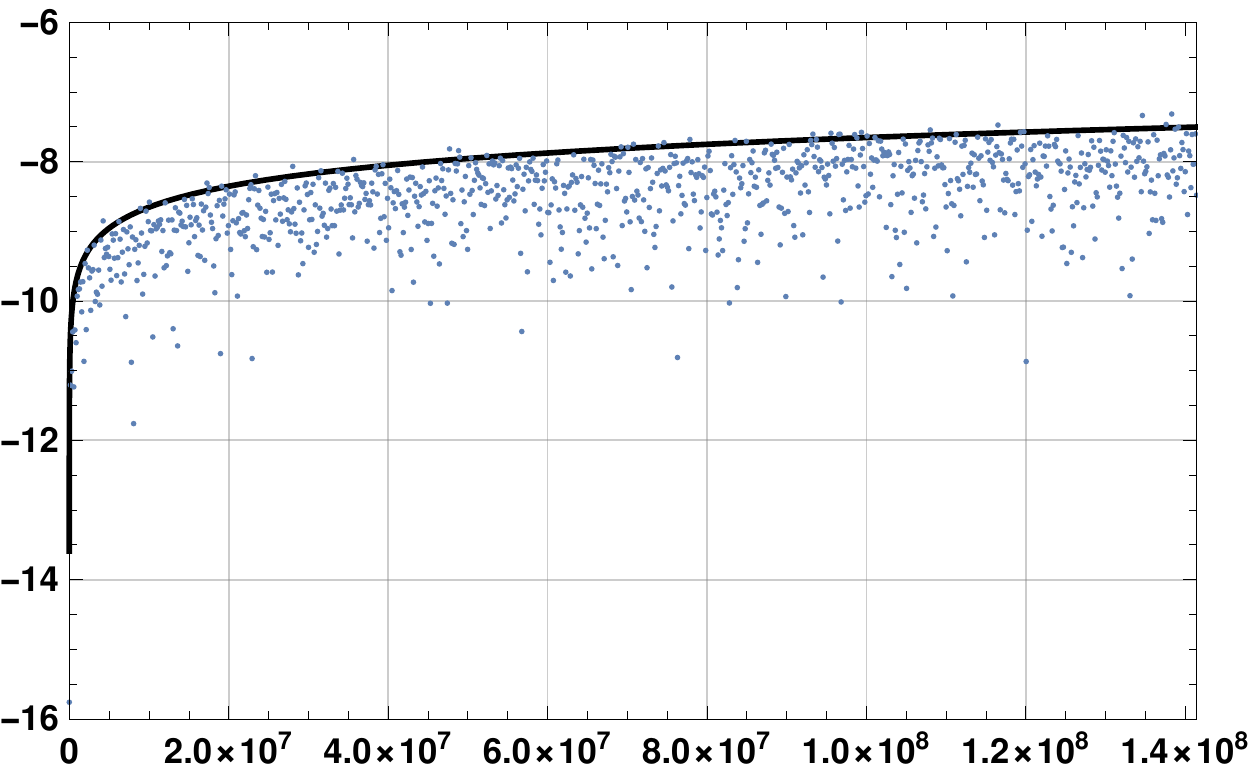}
{\center\small $\nu = 1000\sqrt{2}$}

\end{center}
\caption{The results of the experiments of Section~\ref{experiments:4}.
In each graph, 
the base-$10$ logarithm of the relative errors in calculated values of $H_\nu(t)$ are
plotted as dots 
and the  graph of the function $\log_{10}(\kappa(t)) \epsilon_0$, where
$\kappa(t)$ is the condition number of the evaluation of the function
$H_\nu(t)$ and $\epsilon_0$ is machine epsilon, is plotted as a solid line.
The plots on the left show the results obtained using the {\tt bessel\_eval}
routine while those on the right show results obtained from Amos'
well-known and widely used code \cite{Amos}.
}
\label{figure1}
\end{figure}

\label{experiments:4}
\end{subsection}

\begin{subsection}{The speed and accuracy of the evaluation of $J_n$ and $Y_n$ as a function of $n$}

In these experiments, we compared the speed and accuracy with which
{\tt bessel\_eval} calculates Hankel functions of integer orders
with the speed and accuracy of Amos' code \cite{Amos}.
Reference values were calculated using the well-known three-term recurrence
relations satisfied by Bessel functions.  The calcuation of reference
values was performed using IEEE quadruple precision arithmetic
in order to ensure high accuracy.

\begin{table}[h!!]
\small\center
 \begin{tabular}{l@{\hspace{2em}}cc@{\hspace{2em}}cc}
 \toprule
& \multicolumn{2}{c}{\large \tt bessel\_eval} & \multicolumn{2}{c}{\large \tt Amos' code} \\
\addlinespace[.25em]
$n$ & Maximum relative & Average evaluation   & Maximum relative & Average evaluation \\ 
      & error in $H_n$ & time (in seconds)  & error in $H_n$ & time (in seconds) \\
 \midrule
0 & 7.31\e{-14} & 6.32\e{-07} & 7.20\e{-15} & 4.19\e{-07}  \\
1 & 6.05\e{-13} & 3.24\e{-07} & 5.58\e{-14} & 6.56\e{-07}  \\
10 & 4.10\e{-12} & 3.32\e{-07} & 9.04\e{-13} & 1.07\e{-06}  \\
100 & 4.80\e{-11} & 2.93\e{-07} & 2.13\e{-11} & 2.32\e{-06}  \\
1\sep,000 & 4.51\e{-10} & 2.75\e{-07} & 1.57\e{-10} & 2.12\e{-06}  \\
10\sep,000 & 4.63\e{-09} & 2.73\e{-07} & 3.32\e{-09} & 1.84\e{-06}  \\
100\sep,000 & 4.32\e{-08} & 2.83\e{-07} & 3.08\e{-08} & 1.82\e{-06}  \\
1\sep,000\sep,000 & 4.33\e{-07} & 2.99\e{-07} & - & -  \\
10\sep,000\sep,000 & 4.06\e{-06} & 2.54\e{-07} & - & -  \\
100\sep,000\sep,000 & 2.86\e{-05} & 2.49\e{-07} & - & -  \\
1\sep,000\sep,000\sep,000 & 3.15\e{-04} & 2.39\e{-07} & - & -  \\
 \bottomrule
 \end{tabular}

\caption{The results of the experiments of Section~\ref{experiments:5}
in which
the speed and accuracy with which {\tt bessel\_eval} and the  well-known
code of Amos \cite{Amos} evaluates Hankel functions of integer orders is compared.
Experiments in which Amos' code returned an error code are marked
with dashes.}
\label{table5}
\end{table}

\label{experiments:5}
\end{subsection}

In the first experiment, $n$ was taken to be $0$ and $100$ random points
at which to evaluate $H_n$ were chosen  in the interval $(0,1000).$ 
  In each subsequent 
experiment, $n$ was taken to be a positive integer
and $100$ random points at which to evaluate $H_n$
were chosen from the interval
\begin{equation}
\left(a_n, 1000 n \right),
\end{equation}
where $a_n < \sqrt{n^2-1/4}$ is the solution of  the equation  $\log(-Y_n(a)) = 100$.
In this way, we avoided problems with numerical overflow and underflow.
At each point chosen in this fashion, the 
value of the Hankel function $H_n$ was 
calculated using {\tt bessel\_eval} and with Amos' code.
Table~\ref{table5} reports the results.  There, the maximum observed
relative error in the values of $H_n$ generated by each code is reported
as a function of $n$ as is the average time taken by each code to perform an evaluation.
Amos' code aborts and returns an error code in cases in which it is unable to evaluate
the Bessel functions to at least 7-digit accuracy.  
The corresponding entries  of Table~\ref{table5}
are marked  with dashes.

\begin{subsection}{Extended precision experiments}

It is a straightforward to increase the accuracy of the precomputed
expansions used by the algorithm of this paper.  We constructed a second set
of these expansions, this time asking for 
 $25$ digits of accuracy.  Of course, these precomputations
were conducted using IEEE quadruple precision arithmetic.
We then reran the experiments of Sections \ref{experiments:1}, \ref{experiments:2},
\ref{experiments:3} and  \ref{experiments:5} using IEEE quadruple precision arithmetic instead of the standard
IEEE double precision arithmetic.  Because the laptop we used for experiments
does not support quadruple precision arithmetic in hardware,
it was emulated with software.  This is, of course, highly inefficient and the running times 
of these experiments reflect this fact.  The results 
are shown in Tables~\ref{table1_16} through \ref{table5_16}.

\begin{table}[h!!!]
\begin{center}
\small
 \begin{tabular}{lcc}
 \toprule
 Range of $\nu$ & Maximum relative         & Average evaluation\\
                & error in $\alpha_\nu'(t)$   & time (in seconds) \\
 \midrule
0 - 1 & 2.05\e{-28} & 1.61\e{-05}  \\
1 - 10 & 3.47\e{-28} & 5.86\e{-06}  \\
10 - 100 & 9.62\e{-35} & 5.46\e{-06}  \\
100 - 1\sep,000 & 9.02\e{-29} & 5.45\e{-06}  \\
1\sep,000 - 10\sep,000 & 1.71\e{-32} & 5.60\e{-06}  \\
10\sep,000 - 100\sep,000 & 7.60\e{-33} & 5.43\e{-06}  \\
100\sep,000 - 1\sep,000\sep,000 & 9.62\e{-34} & 5.30\e{-06}  \\
1\sep,000\sep,000 - 10\sep,000\sep,000 & 9.62\e{-35} & 5.29\e{-06}  \\
10\sep,000\sep,000 - 100\sep,000\sep,000 & 9.62\e{-35} & 5.32\e{-06}  \\
100\sep,000\sep,000 - 1\sep,000\sep,000\sep,000 & 2.48\e{-28} & 5.37\e{-06}  \\
 \bottomrule
 \end{tabular}

\end{center}
 \caption{The results obtained by rerunning the  experiments of Section~\ref{experiments:1}
using IEEE quadruple precision arithmetic.  These experiments measure
 the accuracy of the evaluation of $\alpha_\nu'$  in the oscillatory region.
}
 \label{table1_16}
\end{table}

\begin{table}[h!!!]
\small\center
 \begin{tabular}{lccc}
 \toprule
 Range of $\nu$ & Maximum relative    & Maximum relative & Average evaluation\\
                & error in            & error in & time (in seconds) \\
                &$-\nu + \log(J_\nu(t))$ & $\nu + \log(-Y_\nu(t))$   & \\
 \midrule
0.5
 - 1 & 3.83\e{-34} & 2.54\e{-28} & 4.94\e{-05}  \\
1 - 10 & 2.26\e{-25} & 1.18\e{-27} & 1.48\e{-04}  \\
10 - 100 & 7.63\e{-25} & 1.09\e{-27} & 1.10\e{-04}  \\
100 - 1\sep,000 & 8.19\e{-28} & 1.61\e{-27} & 9.87\e{-05}  \\
1\sep,000 - 10\sep,000 & 4.57\e{-28} & 7.80\e{-28} & 5.01\e{-05}  \\
 \bottomrule
 \end{tabular}

 \caption{The results of rerunning the experiments of Section~\ref{experiments:2} using
IEEE quadruple precision arithmetic.  These experiments test
 the accuracy of {\tt bessel\_eval} in the nonoscillatory region
 via comparison with highly accurate reference values
generated using Wolfram's Mathematica package.}
\label{table2_16}
\end{table}

 \begin{table}[h!!]
\center
 \begin{tabular}{lccc}
 \toprule
 Range of $\nu$ & Maximum relative    & Maximum relative & Average evaluation\\
                & error in            & error in & time (in seconds) \\
                &$-\nu + \log(J_\nu(t))$ & $\nu + \log(-Y_\nu(t))$   & \\
 \midrule
100 - 1\sep,000 & 3.00\e{-28} & 2.58\e{-28} & 3.67\e{-05}  \\
1\sep,000 - 10\sep,000 & 2.95\e{-28} & 2.55\e{-28} & 2.77\e{-05}  \\
10\sep,000 - 100\sep,000 & 2.47\e{-28} & 2.38\e{-28} & 2.33\e{-05}  \\
100\sep,000 - 1\sep,000\sep,000 & 2.67\e{-28} & 2.27\e{-28} & 2.02\e{-05}  \\
1\sep,000\sep,000 - 10\sep,000\sep,000 & 2.22\e{-28} & 2.43\e{-28} & 1.90\e{-05}  \\
10\sep,000\sep,000 - 100\sep,000\sep,000 & 1.62\e{-28} & 1.69\e{-28} & 1.74\e{-05}  \\
100\sep,000\sep,000 - 1\sep,000\sep,000\sep,000 & 2.66\e{-28} & 2.76\e{-28} & 1.60\e{-05}  \\
 \bottomrule
 \end{tabular}

 \caption{The results of rerunning the  experiments of Section~\ref{experiments:3} 
using IEEE quadruple precision arithmetic.  These experiments test
 the accuracy with which $-\nu + \log(J_\nu(t))$ and  $\nu + \log(-Y_\nu(t))$  
is evaluated for values of $100 \leq \nu  \leq 1\sep,000\sep,000\sep,000$
and $t \ll \nu$ through comparison with Debye's expansions.}
\label{table3_16}
\end{table}

\begin{table}[h!!]
\begin{center}
\small
 \begin{tabular}{l@{\hspace{2em}}cc@{\hspace{2em}}cc}
 \toprule
& \multicolumn{2}{c}{\large \tt bessel\_eval} & \multicolumn{2}{c}{\large \tt Amos' code} \\
\addlinespace[.25em]
$n$ & Maximum relative & Average evaluation   & Maximum relative & Average evaluation \\ 
      & error in $H_n$ & time (in seconds)  & error in $H_n$ & time (in seconds) \\
 \midrule
0 & 3.00\e{-26} & 5.42\e{-05} & 3.00\e{-17} & 2.53\e{-05}  \\
1 & 2.17\e{-25} & 3.12\e{-05} & 2.85\e{-17} & 1.94\e{-05}  \\
10 & 1.85\e{-24} & 3.23\e{-05} & 1.63\e{-18} & 3.57\e{-05}  \\
100 & 2.39\e{-23} & 2.51\e{-05} & 3.51\e{-18} & 1.96\e{-04}  \\
1\sep,000 & 2.45\e{-22} & 2.14\e{-05} & 3.44\e{-18} & 1.71\e{-04}  \\
10\sep,000 & 8.01\e{-22} & 1.89\e{-05} & 4.64\e{-18} & 1.59\e{-04}  \\
100\sep,000 & 1.48\e{-20} & 1.77\e{-05} & 5.60\e{-18} & 1.59\e{-04}  \\
1\sep,000\sep,000 & 6.08\e{-20} & 1.79\e{-05} & - & -  \\
10\sep,000\sep,000 & 8.52\e{-19} & 2.01\e{-05} & - & -  \\
100\sep,000\sep,000 & 7.62\e{-18} & 1.50\e{-05} & - & -  \\
1\sep,000\sep,000\sep,000 & 5.57\e{-17} & 1.56\e{-05} & - & -  \\
 \bottomrule
 \end{tabular}

\caption{The results of rerunning the experiments of Section~\ref{experiments:6}
using IEEE quadruple precision arithmetic.
In these experments, the speed and accuracy with which {\tt bessel\_eval} and the  well-known
code of Amos \cite{Amos} evaluates Hankel functions of integer orders is compared.
Experiments in which Amos' code returned an error code are marked with dashes.}
\label{table5_16}
\end{center}
\end{table}

\label{experiments:6}
\end{subsection}

\label{section:experiments}
\end{section}

\begin{section}{Conclusions and future work}

Using a simple-minded procedure which can be applied to a large class
of special functions with little modification, we constructed
table which allow for the numerical evaluation of Bessel functions
of nonnegative real orders and positive arguments.  The performance of
the resulting code, at least when it comes to real arguments,
is comparable to that of the well-known and widely used code of 
Amos \cite{Amos}.  

In the nonoscillatory region, our algorithm calculates the logarithms  of the Bessel functions as 
well as their values.  This is useful in cases in which the magnitudes
of the Bessel functions themselves are too large or too small to be encoded
using the IEEE double precision format.
In the oscillatory region, in addition to the values
of the Bessel function itself, our algorithm also returns the values
of a nonoscillatory phase function for Bessel's equation and its derivative.
This is extremely helpful when  computing the zeros of special functions
\cite{BremerZeros}, and when applying special function transforms
via the butterfly algorithm (see, for instance, 
\cite{Candes-Demanet-Ying1,butterfly2,butterfly1,Candes-Demanet-Ying2,Michielssen,ONeil-Rokhlin}).

Preliminary numerical experiments conducted by the author suggest that
adaptive cross approximation (as described in \cite{Bebendorf})
can reduce the size of the table used to evaluate the Bessel functions
by about a factor of $2$  at the cost of a commensurate increase in evaluation time.  Such a trade-off 
is unappealing in the case of Bessel functions given the small size of the precomputed table
needed to evaluate them; however, adaptive cross approximation and related
techniques  might be of great use in other cases, particularly when the special
functions under consideration depend on more than one parameter.

There is a more promising avenue of investigation for reducing the size of the precomputed
expansions.
Well-known results regarding the behavior of solutions of second order differential equations near
their turning points imply that there is a nonoscillatory phase function
which behaves roughly like an error function there.
The portion of the precomputed table dedicated to representing solutions
in this regime is large, and by taking the behavior of solutions near the turning point,
it might be reduced substantially.  Such an approach is being
actively investigated by the author.

The author will also report  on the use of the method of this paper
to evaluate associated Legendre functions and prolate spheroidal wave functions
at a later date, as well as on the rapid application of special function transforms
using techniques related to those discussed here.

\label{section:conclusion}
\end{section}

\begin{section}{References}
\bibliographystyle{acm}
\bibliography{bessel2}
\end{section}
\vfil\eject

\end{document}